
\input amstex
\documentstyle{amsppt}
\loadbold

\magnification=1130
\hsize=6.7truein
\vsize=9.5truein
\hcorrection{-0.1truein}
\vcorrection{-0.2truein}

\document

\baselineskip=14pt

\font\smallbf=cmbx10 at 9pt
\font\smallrm=cmr10 at 9pt
\font\smallit=cmti10 at 9pt
\font\smallsl=cmsl10 at 9pt
\font\sc=cmcsc10

\def \loongrightarrow {\relbar\joinrel\relbar\joinrel\rightarrow}
\def \llongrightarrow
{\relbar\joinrel\relbar\joinrel\relbar\joinrel\rightarrow}
\def \longtwoheadrightarrow
{\relbar\joinrel\relbar\joinrel\twoheadrightarrow}
\def \llongtwoheadrightarrow
{\relbar\joinrel\relbar\joinrel\relbar\joinrel\twoheadrightarrow}
\def \gerg {\frak g}

\def \gerh {\frak h}
\def \gerk {\frak k}
\def \germ {\frak m}

\def \und1 {\underline{1}}

\def \u {{\text{\bf u}}}

\def \h {\hbar}

\def \id {\text{\rm id}}

\def \U {\Cal{U}}
\def \calM {\Cal{M}}

\def \HA {\Cal{HA}}

\def \H {\Bbb{H}}
\def \N {\Bbb{N}}
\def \Z {\Bbb{Z}}

\def \Rhat {\widehat{R}}

\def \Char {\hbox{\it Char}\,}
\def \QrUEA {\Cal{Q}{\hskip1pt}r{\hskip1pt}\Cal{UE{\hskip-1,1pt}A}}
\def \QFA {\Cal{QF{\hskip-1,7pt}A}}


\topmatter

{\ }

\vskip-45pt

   \hfill   {\smallrm {\smallsl Journal f{\"u}r die reine und angewandte Mathematik\/}  (Crelle's Journal)  {\smallbf 612}  (2007), 17--33.}  {\ }
 \vskip1pt
   \hfill   {\smallrm ISSN (Online) 1435-5345, ISSN (Print) 0075-4102, DOI: 10.1515/CRELLE.2007.082, November 2007.}  {\ }

\vskip35pt

\title
  The global quantum duality principle
\endtitle

\vskip-15pt

\author
       Fabio Gavarini
\endauthor


\affil
  Universit\`a degli Studi di Roma ``Tor Vergata'' ---
Dipartimento di Matematica  \\
  Via della Ricerca Scientifica 1, I-00133 Roma --- ITALY  \\
\endaffil

\address\hskip-\parindent
   Fabio Gavarini  \newline
   \indent   Universit{\`a} degli Studi di Roma ``Tor Vergata''  \newline
   \indent   Dipartimento di Matematica  \newline
   \indent   Via della Ricerca Scientifica 1, I-00133 Roma,
ITALY   \newline
   \indent   e-mail: gavarini\@{}mat.uniroma2.it
\endaddress

\abstract
   Let  $ R $  be an integral domain, let  $ \, \h \in R \setminus
\{0\} \, $  be such that  $ \, \Bbbk := R \big/ \h \, R \, $  is
a field, and let  $ \HA $  be the category of torsionless (or
flat) Hopf algebras over  $ R \, $.  We call  $ \, H \in \HA \, $  a
{\it ``quantized function algebra'' (=QFA)}, resp.~{\it ``quantized
restricted universal enveloping algebras'' (=QrUEA)},  at  $ \h \, $
if   --- roughly speaking ---   $ \, H \big/ \h \, H \, $  is the
function algebra of a connected Poisson group, resp.~the (restricted,
if  $ \, R \big/ \h \, R \, $  has positive characteristic) universal
enveloping algebra of a (restricted) Lie bialgebra.  Extending a result
of Drinfeld, we establish an ``inner'' Galois' correspondence on  $ \HA
\, $,  via two endofunctors,  $ (\ )^\vee $  and  $ (\ )' $,  of  $ \HA $
such that  $ H^\vee $  is a QrUEA and  $ H' $  is a QFA  (for all  $ H \!
\in \HA \, $).  In addition:  {\it (a)} \, the image of  $ (\ )^\vee $,
resp.~of  $ (\ )' $,  is the full subcategory of all QrUEAs, resp.~of
all QFAs;  {\it (b)} \,  if  $ \, p := \text{\it Char}(\Bbbk) = 0 \, $,
the restrictions  $ (\ )^\vee{\big|}_{\text{QFAs}} $  and  $ (\ )'
{\big|}_{\text{QrUEAs}} $  yield equivalences inverse to each
other;  {\it (c)} \, if  $ \, p = 0 \, $, starting from a QFA
over a Poisson group  $ G $,  resp.~from a QrUEA over a Lie
bialgebra  $ \gerg \, $,  the functor  $ (\ )^\vee $,
resp.~$ (\ )' $,  gives a QrUEA, resp.~a QFA, over the
dual Lie bialgebra, resp.~the dual Poisson group.  Several,
far-reaching applications are developed in detail in [Ga2--4].
\endabstract

\endtopmatter

\footnote""{Keywords: \ {\sl Hopf algebras, Quantum Groups}.}

\footnote""{ 2000 {\it Mathematics Subject Classification:} \
Primary 16W30, 17B37, 20G42; Secondary 81R50. }

\vskip-21pt

\hfill  \hbox{\vbox{ \baselineskip=10pt
                     \hbox{\smallit  \   ``Dualitas dualitatum }
                     \hbox{\smallit \ \ \;\, et omnia dualitas'' }
                     \vskip4pt
                     \hbox{\smallsl    N.~Barbecue, ``Scholia'' } }
\hskip1,7truecm }

\vskip6pt

\centerline {\bf Introduction }

\vskip10pt

   Generalized ``symmetries'' in mathematics are described by
Hopf algebras.  Among these, the ``geometrical'' ones are of
type  $ \, H = F[G] \, $,  \, the algebra of regular functions
over an algebraic group  $ G $,  and  $ \, H = U(\gerg) \,
\big(\! = \u(\gerg) \big) \, $,  \, the (restricted, if the
ground field  $ \Bbbk $  has positive characteristic) universal
enveloping algebra of a (restricted) Lie algebra  $ \gerg \, $.
These notions of ``geometrical symmetries'' are generalized
by quantum groups: roughly, these are Hopf algebras  $ H $
depending on a parameter  $ \h $  such that, setting  $ \h = 0
\, $,  \, the Hopf algebra one gets is either of the type  $ F[G] $
--- hence  $ H $  is a  {\sl quantized function algebra}, in short
QFA ---   or of the type  $ U(\gerg) $  or  $ \u(\gerg) $  (according
to the characteristic of  $ \Bbbk $)   --- hence  $ H $  is a  {\sl
quantized restricted universal enveloping algebra}, in short QrUEA.
When a QFA exists whose specialization at  $ \, \h = 0 \, $  is
$ F[G] $,  the algebraic group  $ G $  inherits a structure of
Poisson (algebraic) group.  Similarly, if a QrUEA exists whose
specialization is  $ U(\gerg) $  or  $ \u(\gerg) $,  the (restricted)
Lie algebra  $ \gerg $  inherits a structure of Lie bialgebra.  Then,
by general Poisson group theory, Poisson groups $ G^* $  dual to
$ G $  and a Lie bialgebra  $ \gerg^* $  dual to  $ \gerg $  exist.
                                            \par
   In this setting, three basic questions rise at once:
                                            \par
   {\bf --- (1)}  {\sl How can we produce quantum groups?}
                                            \par
   {\bf --- (2)}  {\sl How can we characterize quantum groups
(of either kind) among Hopf algebras?}
                                            \par
   {\bf --- (3)}  {\sl What kind of relationship, if any, does
exist between quantum groups over mutually dual Poisson groups,
or mutually dual Lie bialgebras?}

\vskip1pt

   A first answer to  {\bf (1)}  and  {\bf (3)}  is given, for
$ \, \Char(\Bbbk) = 0 \, $,  by the ``quantum duality principle'',
formulated by Drinfeld in terms of  {\sl formal\/} quantum groups
(cf.~[Dr], \S 7, and [Ga1]): it is a functorial recipe to get, out
of a QFA over  $ G $,  a QrUEA over  $ \gerg^* $,  and a QFA over
$ G^* $  out of a QrUEA over  $ \gerg \, $.
                                            \par
   In this paper I provide a  {\sl global\/}  version of this
principle, which answers questions  {\bf (1)}  through  {\bf (3)}.
Indeed, I push Drinfeld's original method as far as possible, so to apply
it to the category  $ \HA $  of Hopf algebras which are torsion-free (or
flat) over some integral domain, say  $ R \, $,  and to do it for each
$ \, \h \in R \setminus \{0\} \, $  such that  $ \, \Bbbk := R \big/ \h
R \, $  is a field.  In fact, I extend Drinfeld's recipe so to define
endofunctors of  $ \HA \, $.  The image of either functor is contained
in a category of quantum groups (one gives QFAs, the other QrUEAs) so we
answer question  {\bf (1)}.  If  $ \Bbbk $  has zero characteristic, when
restricted to quantum groups these functors yield equivalences inverse to
each other.  Moreover, these equivalences exchange the types of quantum
group (switching QFA with QrUEA) and the underlying Poisson symmetries
(interchanging  $ G $  or  $ \gerg $  with  $ G^* $  or  $ \gerg^* $),
thus solving  {\bf (3)}.  Other details show that these functors endow
$ \HA $  with a (inner) Galois' correspondence, in which QFAs on one
side and QrUEAs on the other side are the subcategories (in  $ \HA $)
of ``fixed points'' for the composition of both Drinfeld's functors (in
suitable order): in particular, this answers question  {\bf (2)}.  Let
me point out that, as my ``Drinfeld's functors'' are defined for each
element  $ \, \h \in R \, $  as above, for any such  $ \h $  and for
any  $ H $  in  $ \HA $  they yield two quantum groups, a QFA and a
QrUEA,  w.r.t.~$ \h $  itself.  Thus we have a method to get, out of
any single  $ H \hskip-2pt \in \hskip-0,3pt \HA \, $, several quantum
groups.
                                                  \par
   Further aspects, examples and applications of the main result are
presented in [Ga2--4].

\vskip4pt

\centerline{ \sc acknowledgements }
  The author thanks P.~Baumann, G.~Carnovale, N.~Ciccoli,
A.~D'Andrea, I.~Damiani, B.~Di Blasio, D.~Fiorenza, L.~Foissy,
A.~Frabetti, C.~Gasbarri and E.~Taft for many helpful discussions.
                                             \par
   A special thank also to the referee for his valuable, fruitful comments and remarks.

\vskip1,3truecm

\centerline {\bf \S \; 1 \ Notation and terminology }

\vskip10pt

  {\bf 1.1 The classical setting.} \, Let  $ \Bbbk $  be a fixed
field of any characteristic.  We call ``algebraic group'' the
maximal spectrum  $ G $  associated to any commutative Hopf
$ \Bbbk $--algebra  $ H \, $; \, then  $ H $  is called the
algebra of regular functions on  $ G $,  denoted  $ F[G] \, $.
We say that  $ G $  is connected if  $ F[G] $  has no non-trivial
idempotents.
%
%
    We denote by  $ \germ_e $  the defining ideal of the unit element
$ \, e \in G \, $  (it is the augmentation ideal of  $ F[G] \, $);
the cotangent space of  $ G $  at  $ e $  is  $ \, \gerg^\times
:= \germ_e \Big/ {\germ_e}^{\!2} \, $,  \,
%
%
which is naturally a Lie coalgebra.  The tangent space of  $ G $
at  $ e $  is the dual space  $ \, \gerg := {\big( \gerg^\times
\big)}^* \, $  to  $ \gerg^\times \, $,  \, which is a Lie algebra.
By  $ U(\gerg) $  we mean the universal enveloping algebra of
$ \gerg \, $:  it is a connected cocommutative Hopf algebra, and
there is a natural Hopf pairing (see  \S 1.2{\it (a)\/})  between
$ F[G] $  and  $ U(\gerg) $.  If  $ \, \hbox{\it Char}\,(\Bbbk) =
p > 0 \, $,  \, then  $ \gerg $  is a restricted Lie algebra, and
$ \, \u(\gerg) := U(\gerg) \Big/ \big( \big\{\, x^p - x^{[p\hskip0,7pt]}
\,\}_{x \in \gerg} \big) \, $  is the restricted universal enveloping
algebra of  $ \gerg \, $.  To unify notation and terminology, when
$ \, \Char(\Bbbk) = 0 \, $  we call any Lie algebra  $ \gerg $
``restricted'', by its ``restricted universal enveloping algebra''
we mean  $ U(\gerg) $,  and we write  $ \, \U(\gerg) := U(\gerg) \, $
if  $ \, \Char(\Bbbk) = 0 \, $  and  $ \, \U(\gerg) := \u(\gerg) \, $
if  $ \, \Char(\Bbbk) > 0 \, $.
                                            \par
   Let  $ H $  be a Hopf algebra over an integral domain  $ D \, $.
We call  $ H $  a ``function algebra'' (FA in short) if it is commutative,
with no non-trivial idempotents, and such that, if  $ \, p := \Char(\Bbbk)
> 0 \, $,  \, then  $ \, \eta^p = 0 \, $  for all  $ \eta $  in the kernel
of the counit of  $ H \, $.  If  $ D $  is a field, an FA is the algebra
of regular functions of an algebraic group-scheme over  $ D $  which is
connected and, if  $ \, \Char(\Bbbk) > 0 \, $,  \, is zero-dimensional
of height 1; conversely, if  $ G $  is such a group-scheme then  $ F[G] $
has these properties.  Instead, we call  $ H $  a ``restricted universal
enveloping algebra'' (=rUEA) if it is cocommutative, connected, and
generated by its primitive part.  If  $ D $  is a field, an rUEA is
the restricted universal enveloping algebra of some (restricted) Lie
algebra over  $ D \, $;  \, conversely, if  $ \gerg $  is such a Lie
algebra, then  $ \U(\gerg) $  has these properties
%
%
(see, e.g., [Mo], Theorem 5.6.5, and references therein).
                                            \par
   Now assume  $ G $  is a Poisson group (for this and other notions
hereafter see [CP], but within an  {\sl algebraic geometry\/}  setting).
Then  $ F[G] $  is a Poisson Hopf algebra, whose Poisson bracket
induces on  $ \gerg^\times $  a Lie bracket which makes it into a Lie
bialgebra; hence  $ U(\gerg^\times) $  is a co-Poisson Hopf algebra too.
Also,  $ \gerg $  is a Lie bialgebra (in topological sense, if  $ G $
is infinite dimensional) too, and  $ U(\gerg) $  is a (maybe topological)
co-Poisson Hopf algebra.  The Hopf pairing between  $ F[G] $  and
$ U(\gerg) $  then is compatible with these additional co-Poisson and
Poisson structures.  Moreover, the perfect (=non-degenerate) evaluation
pairing  $ \, \gerg \times \gerg^\times \! \rightarrow \Bbbk \, $
is compatible with the Lie bialgebra structure on either side (see \S
1.2{\it (b)\/}):  so  $ \gerg $  and  $ \gerg^\times $  are Lie bialgebras  {\sl dual to each other}.  In the sequel, we denote by
$ G^\star $  any connected algebraic Poisson group with  $ \gerg $  as cotangent Lie bialgebra, and say
   \hbox{it is  {\sl dual\/}  to  $ G \, $.}
                                            \par
   If  $ H $  is a Hopf algebra we denote its Hopf operations
by  $ \Delta $  (the coproduct),  $ \epsilon $  (the counit) and
$ S $  (the antipode),  and we use standard  $ \sigma $--notation
$ \, \Delta(x) = \sum_{(x)} x_{(1)} \otimes x_{(2)} \, $  for all
$ \, x \in H \, $.

\vskip9pt

\proclaim {Definition 1.2}
                                   \hfill\break
   \indent   (a) \, Let  $ H $,  $ K $  be Hopf algebras (maybe in
topological sense) over a ring  $ R \, $.  A  {\sl Hopf (algebra)
pairing\/}  between  $ H $  and  $ K $  is a pairing  $ \; \langle
\,\ , \,\ \rangle \, \colon \, H \times K \! \longrightarrow R \; $
such that  $ \;\; \big\langle x, y_1 \, y_2 \big\rangle = \big\langle
\Delta(x), y_1 \otimes y_2 \big\rangle := \sum_{(x)} \big\langle x_{(1)},
y_1 \big\rangle \, \big\langle x_{(2)}, y_2 \big\rangle \, $,  $ \;
\big\langle x_1 \, x_2, y \big\rangle = \big\langle x_1 \otimes x_2,
\Delta(y) \big\rangle := \sum_{(y)} \big\langle x_1, y_{(1)} \big\rangle
\, \big\langle x_2, y_{(2)} \big\rangle  \, $,  $ \, \langle x, 1 \rangle
= \epsilon(x) \, $,  $ \; \langle 1, y \rangle = \epsilon(y) \, $,  $ \;
\big\langle S(x), y \big\rangle = \big\langle x, S(y) \big\rangle \, $,
\, for all  $ \, x, x_1, x_2 \in H $,  $ \, y, y_1, y_2 \in K $.
                                   \hfill\break
   \indent   (b) \, Let  $ \gerg $,  $ \gerh $  be Lie bialgebras (maybe
in topological sense) over a ring  $ \Bbbk \, $.  A  {\sl Lie bialgebra
pairing\/}  between  $ \gerg $  and  $ \gerh $  is a pairing  $ \; \langle
\,\ , \,\ \rangle \, \colon \, \gerg \times \gerh \longrightarrow \Bbbk
\; $  such that  $ \;\; \big\langle x, [y_1,y_2] \big\rangle = \big\langle
\delta(x), y_1 \otimes y_2 \big\rangle := $\break
 $ \sum_{[x]} \! \big\langle x_{[1]},
y_1 \big\rangle \, \big\langle x_{[2]}, y_2 \big\rangle  \, $,  $ \;
\big\langle [x_1,x_2], y \big\rangle = \big\langle x_1 \otimes x_2,
\delta(y) \big\rangle := \sum_{[y]} \! \big\langle x_1, y_{[1]} \big\rangle
\, \big\langle x_2, y_{[2]} \big\rangle \, $,  \, for all  $ \, x, x_1,
x_2 \! \in \gerg \, $  and  $ \, y, y_1, y_2 \in \gerh $,  \, with
$ \, \delta(x) = \sum_{[x]} x_{[1]} \otimes x_{[2]} \, $  and
$ \, \delta(y) = \sum_{[y]} y_{[1]} \otimes y_{[2]} \, $.
\endproclaim

\vskip5pt

  {\bf 1.3 The quantum setting.} \, Let  $ R $  be an integral domain,
$ \, F = F(R) \, $  its field of quotients.  Let  $ \calM $  be the category of torsion-free  $ R $--modules,  $ \HA $  the category of all Hopf algebras in  $ \calM \, $.  Let  $ \calM_F $  be the category of  $ F $--vector  spaces,  $ \HA_F $  the category of all Hopf algebras in  $ \calM_F \, $;  \, for  $ \, M \in \calM \, $,  \, set  $ \, M_F := F(R) \otimes_R M \, $.  A subset  $ \, H \subseteq \Bbb{H} \in \HA_F \, $  is called  {\it an  $ R $--integer  form}  (or  {\it an  $ R $--form})  {\it of\/}  $ \, \Bbb{H} $  iff  $ \, H \, $  is a Hopf  $ R $--subalgebra  of  $ \, \Bbb{H} \, $  (hence in particular  $ \, H \in \HA \, $)  and  $ \; H_F := F(R) \otimes_R H = \Bbb{H}
\, $.
                                            \par
   Let  $ \, \h \in R \setminus \{0\} \, $  be prime (fixed
throughout), and  $ \, \Bbbk := R \big/ (\h) = R \big/ \h \, R
\, $.  For any  $ R $--module  $ M $,  set  $ \, M_\h{\Big|}_{\h=0}
\!\! := M \big/ \h \, M = \Bbbk \otimes_R M \, $  (the  {\sl
specialization\/}  of  $ M $  at  $ \, \h = 0 \, $)  and  $ \,
M_\infty := \bigcap_{n=0}^{+\infty} \h^n M \, $.  Finally, for
any  $ H \! \in \! \HA \, $,  let  $ \, I_{\scriptscriptstyle H}
:= \text{\sl Ker} \Big( H \,{\buildrel \epsilon \over
{\relbar\joinrel\relbar\joinrel\twoheadrightarrow}}\,
R \, {\buildrel {\h \mapsto 0} \over
{\relbar\joinrel\relbar\joinrel\relbar\joinrel\twoheadrightarrow}}
\, \Bbbk \Big) \, $,  \, and set  $ \, {I_{\scriptscriptstyle
H}}^{\!\infty} := \bigcap_{\,n=0}^{+\infty}
{I_{\scriptscriptstyle H}}^{\!n} \, $.

\vskip9pt

\proclaim {Definition 1.4} ({\sl ``Global quantum groups''})
Let  $ \, \h \in R \setminus \{0\} \, $  be a prime, and
$ \; \Bbbk := R \big/ \h \, R \; $.
                                               \par
   (a) \, We call  {\sl quantized restricted universal
enveloping algebra (at  $ \h $)}   --- in short,  {\sl QrUEA\/}
---   any  $ \U_\h \! \in \! \HA $  such that  $ \, \U_\h{\big|}_{\h=0} \! := \U_\h \big/ \h \, \U_\h $
is a restricted universal
   \hbox{enveloping algebra \! (an rUEA) \! over  $ \Bbbk $.}
                                               \par
   We call  $ \, \QrUEA \, $  the full subcategory of  $ \, \HA $  whose objects are all the QrUEAs (at  $ \h $).
                                               \par
   (b) \, We call  {\sl quantized function algebra (at  $ \h $)}
--- in short,  {\sl QFA\/}  ---   any  $ \, F_\h \in \HA $  such
that  $ \, F_\h{\big|}_{\h=0} \! := F_\h \big/ \h \, F_\h \, $
is a function algebra (an FA) over\/  $ \Bbbk \, $,  \,  and  $ \, {(F_\h)}_\infty = {I_{\!\scriptscriptstyle F_\h}}^{\!\!\infty}
\, $  (notation of \S 1.3).
                                               \par
   We call  $ \, \QFA \, $  the full subcategory of  $ \, \HA $  whose objects are all the QFAs (at  $ \h $).
\endproclaim
 \eject

\vskip1pt

   {\bf Remark 1.5:} \, If  $ \, \U_\h \, $  is a QrUEA (at
$ \h \, $)  then  $ \, \U_\h{\big|}_{\h=0} \, $  is a co-Poisson
Hopf algebra, w.r.t.~the Poisson cobracket  $ \delta $  defined as follows: if  $ \, x \in \U_\h{\big|}_{\h=0} \, $  and  $ \, x' \in \U_\h \, $  gives  $ \, x = x' \mod \h \, \U_\h \, $,  \, then  $ \, \delta(x) := \big( \h^{-1} \, \big( \Delta(x') - \Delta^{\text{op}}
(x') \big) \big) \mod \h \, \big( \U_\h \otimes \U_\h \big) \, $.
So, if  $ \Bbbk $  is a field, then  $ \, \U_\h{\big|}_{\h=0} \cong \U(\gerg) \, $  for some Lie algebra  $ \gerg \, $,  and by [Dr],
\S 3, the restriction of  $ \delta $  makes  $ \gerg $  into a  {\sl Lie bialgebra\/};
%
%
 then I shall write  $ \, \U_\h = \U_\h(\gerg) \, $.  Similarly,
if  $ F_\h $  is a QFA at  $ \h $,  then  $ \, F_\h{\big|}_{\h=0}
\, $  is a  {\sl Poisson\/}  Hopf algebra, w.r.t.~the Poisson
bracket  $ \{\,\ ,\ \} $  defined as follows: if  $ \, x $,
$ y \in F_\h{\big|}_{\h=0} \, $  and  $ \, x' $,  $ y' \in
F_\h \, $  give  $ \, x = x' \mod \h \, F_\h \, $,  $ \, y =
y' \mod \h \, F_\h \, $,  \, then  $ \, \{x,y\} := \big( \h^{-1}
(x' \, y' - y' \, x') \big) \mod \h \, F_\h \, $.  Thus, if
$ \Bbbk $  is a field,  $ \, F_\h{\big|}_{\h=0} \cong F[G] \, $
for some connected  {\sl Poisson\/}  algebraic group  $ G \, $:
%
%
 in this case I shall write
$ \, F_\h = F_\h[G] \, $.

%
%

\vskip6pt

\proclaim{Definition 1.6}
                                           \hfill\break
   \indent   (a) \, Let  $ R $  be an integral domain, and let
$ F $  be its field of fractions.  Given two  $ F $--modules
$ \Bbb{A} $,  $ \Bbb{B} $,  and an  $ F $--bilinear  pairing
$ \; \Bbb{A} \times \Bbb{B} \longrightarrow F \, $,  \; for any
$ R $--submodule  $ \, A \subseteq \Bbb{A} \, $  and  $ \, B
\subseteq \Bbb{B} \, $  we define  $ \; \displaystyle{ A^\bullet
\, := \big\{\, b \in \Bbb{B} \;\big\vert\; \big\langle A, \, b
\big\rangle \subseteq R \big\} } \; $  and  $ \; \displaystyle{
B^\bullet \, := \big\{\, a \in \Bbb{A} \;\big\vert\; \big\langle
a, B \big\rangle \subseteq R \big\} } \, $.
                                           \hfill\break
   \indent   (b) \, Let  $ R $  be an integral domain.  Given
$ \, H $,  $ K \in \HA \, $,  \, we say that  {\sl  $ H $  and
$ K $  are dual to each other}  if there exists a perfect Hopf
pairing between them for which  $ \, H = K^\bullet \, $  and
$ \, K = H^\bullet \; $.
\endproclaim

\vskip1,3truecm

\centerline {\bf \S \; 2 \  The global quantum duality principle }

\vskip13pt

  {\bf 2.1 Drinfeld's functors.} \,  (Cf.~[Dr], \S 7) Let  $ R \, $,
$ F \, $,  $ \HA $  and  $ \, \h \in R \setminus \{0\} \, $  be
as in \S 1.3.  For any  $ \, H \in \HA \, $,  \, let  $ \, I =
I_{\scriptscriptstyle \! H} := \hbox{\sl Ker} \Big( H \,{\buildrel
\epsilon \over {\relbar\joinrel\twoheadrightarrow}}\, R \,{\buildrel
{\h \mapsto 0} \over \llongtwoheadrightarrow}\, R \big/ \h \, R =
\Bbbk \Big) = \hbox{\sl Ker} \Big( H \,{\buildrel {\h \mapsto 0}
\over \llongtwoheadrightarrow}\, H \big/ \h \, H \,{\buildrel
\bar{\epsilon} \over {\relbar\joinrel\twoheadrightarrow}}\,
\Bbbk \Big) \, $,  \, as in \S 1.3, where  $ \bar{\epsilon} $
denotes the counit of  $ \, H{\big|}_{\h=0} \; $.  I define
  $$  H^\vee \; := \; {\textstyle \sum_{n \geq 0}} \, \h^{-n} I^n \; =
\; {\textstyle \sum_{n \geq 0}} {\big( \h^{-1} I \, \big)}^n \; = \;
{\textstyle \bigcup_{n \geq 0}} {\big( \h^{-1} I \, \big)}^n \quad
\big( \! \subseteq H_F \, \big) \; .  $$
If  $ \, J = J_{\scriptscriptstyle \! H} := \hbox{\sl Ker}\,
(\epsilon_{\scriptscriptstyle \! H}) \, $  then  $ \, I = J
+ \h \, R \cdot 1_{\scriptscriptstyle H} \, $,  \, so  $ \;
H^\vee = \sum_{n \geq 0} \h^{-n} J^n = \sum_{n \geq 0}
{\big( \h^{-1} J \, \big)}^n \; $  too.
                                              \par
  Given any Hopf algebra  $ H $,  for every  $ \, n \in \N \, $  define  $ \; \Delta^n \colon H \longrightarrow H^{\otimes n}
\; $  by  $ \, \Delta^0 := \epsilon \, $,  $ \, \Delta^1 := \id_{\scriptscriptstyle H} $,  \, and  $ \, \Delta^n := \big(
\Delta \otimes \id_{\scriptscriptstyle H}^{\otimes (n-2)}
\big) \circ \Delta^{n-1} \, $  if  $ \, n > 2 \, $.  For any
ordered subset  $ \, \Sigma = \{i_1, \dots, i_k\} \subseteq
\{1, \dots, n\} \, $  with  $ \, i_1 < \dots < i_k \, $,  \, define
the morphism  $ \; j_{\scriptscriptstyle \Sigma} : H^{\otimes k}
\longrightarrow H^{\otimes n} \; $  by  $ \; j_{\scriptscriptstyle
\Sigma} (a_1 \otimes \cdots \otimes a_k) := b_1 \otimes \cdots
\otimes b_n \; $  with  $ \, b_i := 1 \, $  if  $ \, i \notin
\Sigma \, $  and  $ \, b_{i_m} := a_m \, $  for  $ \, 1 \leq m
\leq k \, $.  Set  $ \; \Delta_\Sigma := j_{\scriptscriptstyle
\Sigma} \circ \Delta^k \, $,  $ \, \Delta_\emptyset := \Delta^0
\, $,  and  $ \; \delta_\Sigma := \sum_{\Sigma' \subset \Sigma}
{(-1)}^{n- \left| \Sigma' \right|} \Delta_{\Sigma'} \, $,  $ \;
\delta_\emptyset := \epsilon \, $;  \, this admits the inverse formula
$ \; \Delta_\Sigma = \sum_{\Psi \subseteq \Sigma} \delta_\Psi \, $.
We shall use notation  $ \, \delta_0 := \delta_\emptyset \, $,  $ \,
\delta_n := \delta_{\{1, 2, \dots, n\}} \, $,  and the useful formula
$ \, \delta_n = {(\id_{\scriptscriptstyle H} \! - \epsilon)}^{\otimes n}
\! \circ \Delta^n \;\; (n \in \N_+) \, $.
                                               \par
   Now consider again any  $ \, H \in \HA \, $  and  $ \, \h \in
R \setminus \{0\} \, $  as in \S 1.3.  Then I define
  $$  H' \; := \; \big\{\, a \in H \,\big\vert\, \delta_n(a) \in
\h^n H^{\otimes n} , \; \forall \,\, n \in \N \, \big\}  \quad
\big( \! \subseteq H \, \big) \, .  $$
   \indent   Now I can state the main result of the paper:

\vskip11pt

\proclaim {Theorem 2.2} ({\sl ``The Global Quantum Duality Principle''})
Assume that  $ \, \Bbbk := R \big/ \h \, R \, $  is a field.
                                        \hfill\break
  \indent   (a) \, The assignment  $ \, H \mapsto H^\vee \, $,
resp.~$ \, H \mapsto H' \, $,  defines a functor  $ \; {(\ )}^\vee
\colon \, \HA \llongrightarrow \HA \, $,  \, {resp.}  $ \; {(\ )}' \colon
\, \HA \llongrightarrow \HA \, $,  \, whose image lies in  $ \QrUEA $,
resp.~in  $ \QFA $.  Moreover, for all  $ \, H \! \in \HA \, $  we
have  $ \, H \! \subseteq \! {\big( H^\vee \big)}' \, $  and  $ \; H
\! \supseteq \! {\big(H'\big)}^{\!\vee} \, $,  \; hence  $ \; H^\vee
\! = \! \big(\big(H^\vee\big)' \,\big)^{\!\vee} \, $  and  $ \; H' \!
= \! \big( \big(H'\big)^{\!\vee} \big)' \; $.  In addition, if  $ \, H
\in \HA \, $  is flat, then  $ H^\vee $  and  $ H' $  are flat as well.
                                        \hfill\break
   \indent   (b) \, Assume that  $ \; \hbox{\it Char}\,(\Bbbk)
= 0 \, $.  Then for any  $ \, H \in \HA \, $  we have
 \vskip1pt
   \centerline{ $ \displaystyle{ H = {\big(H^\vee\big)}'
\,\Longleftrightarrow\, H \in \QFA  \qquad  \hbox{and}
\qquad  H = {\big( H' \big)}^{\!\vee} \,\Longleftrightarrow\,
H \in \QrUEA } $ }
 \vskip1pt
\noindent   thus  $ {(\ )}^\vee $  and  $ {(\ )}' $  restrict
to equivalences, inverse to each other, between  $ \QFA $  and
$ \QrUEA \, $.
%
%
 \eject
  \indent   (c) \, (``Quantum Duality Principle'') Assume that
$ \; \hbox{\it Char}\,(\Bbbk) = 0 \, $.  Then
 \vskip2pt
   \centerline{ $ {F_\h[G]}^\vee{\Big|}_{\h=0} := {F_\h[G]}^\vee \Big/
\h \, {F_\h[G]}^\vee \, = \, U(\gerg^\times) \; ,  \quad
\displaystyle{ {U_\h(\gerg)}'{\Big|}_{\h=0} := {U_\h(\gerg)}' \Big/
\h \, {U_\h(\gerg)}' \, = \, F\big[G^\star\big] } $ }
 \vskip3pt
\noindent
 (cf.~\S\S 1.1, 1.5)
where the choice of  $ \, G^\star $  (among all the connected Poisson
algebraic groups with tangent Lie bialgebra  $ \gerg^\star $)  depends
on the choice of  $ \, U_\h(\gerg) \, $.  In other words,  $ \,
{F_\h[G]}^\vee \, $  is a QrUEA for the Lie bialgebra
$ \gerg^\times $,  and  $ \, {U_\h(\gerg)}' \, $  is
a QFA for the Poisson group  $ G^\star $.
                                        \hfill\break
  \indent   (d) \, Assume that  $ \, \hbox{\it Char}\,(\Bbbk) = 0 \, $.
Let  $ \, F_\h \in \QFA \, $,  $ \, U_\h \in \QrUEA \, $  be dual to
each other ({w.r.t.} some pairing).  Then  $ \, {F_\h}^{\!\!\vee} $  and
$ \, {U_\h}' $  are dual to each other (w.r.t.~the  {\sl same}  pairing).
                                        \hfill\break
  \indent   (e) \, Assume that  $ \, \hbox{\it Char}\,(\Bbbk) = 0 \, $.
Then for any  $ \, \H \in \HA_F \, $  the following are equivalent:
                                        \hfill\break
   \indent \indent  (e--1) \;  $ \H $  has an  $ R $--integer  form
$ H_{(f)} $  which is a QFA at  $ \h \, $;
                                        \hfill\break
   \indent \indent  (e--2) \;  $ \H $  has an  $ R $--integer  form
$ H_{(u)} $  which is a QrUEA at  $ \h \, $.
\endproclaim

\vskip7pt

{\bf Remarks 2.3:}  {\it (a)} \, In [Ga2] the effect of Drinfeld's
functors on some popular quantum groups or other Hopf algebras is
studied in detail.  An important application to ``classical'' Hopf
algebras is explained in [Ga3]: a special case of it, regarding the
Nottingham group, is studied in [Ga4].
                                                  \par
   {\it (b)} \, Theorem 2.2 can be still partially generalized, see
\S 4.10 at the end of the paper.

 \vskip1,3truecm

\centerline {\bf \S \; 3 \  General properties of Drinfeld's functors }

\vskip13pt

   We begin with a few technicalities, then pass to the first
relevant results.  Fix  $ R $  and  $ \h $  as in \S 1.3.

\vskip9pt

\proclaim{Lemma 3.1} \, Let  $ \, H \in \HA \, $,  \, and set  $ \,
\overline{H} := H \big/ H_\infty \, $  (notation of \S 1.3).  Then:
                                         \hfill\break
   \indent   {\it (a)}  $ \, H_\infty \! = {(H')}_\infty \, $,
$ \, H_\infty \! \subseteq \! {\big(H^\vee\big)}_\infty \, $,
$ \, H_\infty $  is a Hopf ideal
       \hbox{and subcoalgebra of  $ H \, $,  and
$ {\big( \overline{H} \,\big)}_\infty \!\! = \! \{0\} $.}
 Moreover, there are
natural isomorphisms  $ \; {\big( \overline{H} \,\big)}^{\!\vee}
\cong H^\vee \! \Big/ H_\infty \; $  and  $ \; {\big( \overline{H}
\,\big)}' \cong H' \! \Big/ H_\infty \; $.
                                         \hfill\break
   \indent   {\it (b)}  $ \; \overline{H} \in \HA \, $,  \, and
$ \, {\overline{H}}{\big|}_{\h=0} = H{\big|}_{\h=0} \; $.  In
particular, if  $ H $  is a QFA, then  $ \overline{H} $  is
a QFA too, and if  $ H $  is a QrUEA, then  $ \overline{H} $
is a QrUEA too.  \qed
\endproclaim

\vskip7pt

\proclaim{Lemma 3.2} \! ([KT], Lemma 3.2) Let  $ \, H \! \in \!
\HA \, $,  $ \, a $,  $ b \! \in \! H \, $,  \, and  $ \, \Phi
\! \subseteq \! \N \, $  a  {\sl finite}  subset.  Then\break
   \indent  (a) \quad \;\;  $ \displaystyle{
\delta_\Phi(ab) = {\textstyle \sum_{\Lambda \cup Y = \Phi}}
\, \delta_\Lambda(a) \, \delta_Y(b) } \quad $;
                                     \hfill\break
   \indent  (b) \quad  if  $ \, \Phi \not= \emptyset \, $,
\; then  \quad  $ \displaystyle{ \delta_\Phi(ab - ba) =
{\textstyle \sum_{\Sb  \Lambda \cup Y = \Phi  \\
\Lambda \cap Y \not= \emptyset  \endSb}}
\big( \delta_\Lambda(a) \, \delta_Y(b) - \delta_Y(b) \,
\delta_\Lambda(a) \big) } \quad $.   \qed
\endproclaim

\vskip6pt

  {\sl From now on, we make the following\/}  {\bf assumption:
$ \, \Bbbk := R \big/ \h \, R \, $  is a field}.

\vskip10pt

\proclaim{Proposition 3.3}  $ \, H \mapsto \! H^\vee $
and  $ \, H \mapsto \! H' $  gives well-defined endofunctors
of  $ \, \HA \, $,  \, which preserve flatness.  They enjoy
$ \, H \! \subseteq \! {\big(H^\vee\big)}' $,  $ \, H \!
\supseteq \! {\big( H' \big)}^{\!\vee} \! $,  so  $ \,
H^\vee \! = \! \Big( \! {\big(H^\vee\big)}' \Big)^{\!\!\vee} $,
$ \, H' \! = \! \Big( \! {\big(H'\big)}^{\!\vee} \Big)' $,  \,
for all  $ \, H \in \HA \, $.
\endproclaim

\demo{Proof}  Given  $ \, H \in \HA \, $,  \, clearly  $ H^\vee $
and  $ H' $  are torsion-free, hence  $ \, H^\vee, H' \in \calM
\, $.  Moreover, flatness is preserved when taking submodules and/or
localizing, so  $ H^\vee $  and  $ H' $  are flat is  $ H $  is flat.
                                              \par
   Since  $ \, J := \hbox{\sl Ker}\,(\epsilon_{\scriptscriptstyle
\! H}) \, $  is a Hopf ideal of  $ H \, $,  \, we see at once
that  $ \, H^\vee \in \HA \, $  also.
                                              \par
   On the other hand,  $ H' $  is a unital  $ R $--subalgebra  of
$ H $  thanks to Lemma 3.2{\it (a)\/}  and the very definitions.
Moreover,  $ \, \Delta^n \circ S = S^{\otimes n} \circ \Delta^n
\, $  implies  $ \, \delta_n \circ S = S^{\otimes n} \circ \delta_n
\, $  ($ \, n \in \N \, $),  \, which yields  $ \, S(H') = H' \, $.
Thus we still need only to prove that  $ H' $  is a subcoalgebra,
namely  $ \, \Delta\big(H'\big) \subseteq H' \otimes H' \, $.
                                        \par
   Consider  $ \; D_H^{\;\prime} := \big\{\, z \in H \otimes H
\;\big|\; (\delta_r \otimes \delta_s)(z) \in \h^{r+s} H^{\otimes r}
\otimes H^{\otimes s}, \; \forall\, r, s \in \N \,\big\} \, $.  The
coassociativity of  $ \Delta $  yields  $ \, (\delta_r \otimes \delta_s)
\circ \Delta = \delta_{r+s} \, $,  \, whence we get  $ \, \Delta \big(
H' \big) \subseteq D_H^{\;\prime} \, $.  Therefore, the claim will
follow once we prove the identity  $ \; D_H^{\;\prime} = H' \otimes
H' \; $.  The key point is to prove that  $ \; D_H^{\;\prime} \subseteq
H' \otimes H' \; $,  \; the converse inclusion being trivial; to this
end, we shall resort to completions.  Note that, by Lemma 3.1, we can
reduce to prove the main statement for  $ \, \overline{H}^{\,\prime} \, $.
Even more, one clearly has  $ \; D_H^{\;\prime} = H' \otimes H' \; $  if
and only if  $ \; D_{\overline{H}}^{\;\prime} = \overline{H}^{\,\prime}
\otimes \overline{H}^{\,\prime} \; $.  Therefore, we can assume that
$ \, H_\infty = \{0\} \, $.
                                        \par
   Let  $ \, \Rhat \, $  and  $ \, \widehat{H} \, $  respectively be the
$ \h $--adic  completion of  $ R $  and  $ H \, $:  then  $ \, \widehat{H}
\, $  is a separated complete topological  $ \Rhat $--module,  hence it
is topologically free, and is a topological Hopf algebra.  Its coproduct
takes values into the  $ \h $--adic  completion  $ \, H \, \widehat{\otimes}
\, H \, $  of  $ \, H \otimes H \, $,  \, and  $ \, H \, \widehat{\otimes}
\, H = \widehat{H} \, \widehat{\otimes} \, \widehat{H} = \widehat{H
\otimes H} \, $.  As  $ H $  embeds into  $ \widehat{H} $  (because
$ \, H_\infty = \{0\} \, $),  we identify  $ H $  itself with a (dense)
Hopf  $ R $--subalgebra  of  $ \widehat{H} \, $.  Then  $ \; \h^n H = H
\bigcap  \big( \h^n \widehat{H} \,\big) \; $  and  $ \; \widehat{H} \Big/
\h^n \widehat{H} = H \Big/ \h^n H \, $,  \; for all  $ \, n \in \N \, $.
                                        \par
   Let  $ \; \widehat{H}' := \Big\{\, a \in \widehat{H} \;\Big\vert\;
\delta_n(a) \in \h^n \widehat{H}^{\widehat{\otimes} \, n} , \; \forall
\,\, n \in \N \, \Big\} \, $.  Then clearly  $ \, H' = \widehat{H}'
\bigcap H \, $.  Similarly, define  $ \big( \widehat{H \otimes H}
\big)' $  simply taking  $ H \otimes H $  instead of  $ H \, $,  \,
and define  $ \, D_{\widehat{H}}^{\;\prime} \, $  just mimicking the
definition of  $ D_H^{\;\prime} \, $.
                                        \par
   Now, Proposition 2.1 in [EH] proves that  $ \; D_{\widehat{H}}^{\;
\prime} = \big( \widehat{H \otimes H} \big)' \; $;  \; to be precise, the statement in [EH] is for quantized universal enveloping algebras, but the arguments used therein can be easily modified as to apply
to our  $ \widehat{H} $  as well.  In addition, one has also  $ \,
\big( \widehat{H \otimes H} \big)' = \widehat{H}' \,\bar{\otimes}\, \widehat{H}' \, $,  \, where  $ \, \bar{\otimes} \, $  denotes a topological tensor product (w.r.t.~the weak topology, see [EH] or [Ga1] for details), hence one gets  $ \, \big( \widehat{H \otimes H} \big)' \bigcap \; (H \otimes H) = \big( \widehat{H}' \,\bar{\otimes}\, \widehat{H}' \big) \bigcap \; (H \otimes H) = \big( \widehat{H}' \bigcap \, H \big) \otimes \big( \widehat{H}' \bigcap \, H \big) = H' \otimes H' \, $.  But then, noting that  $ \, D_H^{\;\prime} = D_{\widehat{H}}^{\;\prime} \, \bigcap \, \big( H \otimes H \big) \, $,  \, the identity
$ \; D_{\widehat{H}}^{\;\prime} = \big( \widehat{H \otimes H} \big)'
\; $  eventually yields  $ \; D_H^{\;\prime} = D_{\widehat{H}}^{\;\prime} \, \bigcap \, \big( H \otimes H \big) = \big( \widehat{H \otimes H} \big)' \, \bigcap \, \big( H \otimes H \big) = H' \otimes H' \, $,  \; q.e.d.
                                              \par
   The outcome is that  $ \, H' \!\in \HA \, $,  \, hence  $ (\ )' $
is well defined on objects and takes values into  $ \HA \, $.
                                           \par
   As to morphisms, let  $ \, \varphi \in \! \hbox{\sl
Mor}_{\scriptscriptstyle \HA}(H,K) \, $.
%
%
Then its scalar extension  $ \varphi_{F(R)} $  defines  $ \; \varphi^\vee
:= \varphi_{F(R)}{\big|}_{H^\vee} \in \hbox{\sl Mor}_{\scriptscriptstyle
\HA}\big( H^\vee, \, K^\vee \big) \, $;  \, similarly one defines  $ \,
\varphi' := \varphi{\big\vert}_{H'} \in \hbox{\sl Mor}_{\scriptscriptstyle
\HA}\big(H',K'\big) \, $  as well.
                                                 \par
   Finally, let  $ \, H \in \HA \, $.  For  $ \, n \in \N \, $
we have  $ \, \delta_n(H) \subseteq {J_{\scriptscriptstyle H}}^{\!
\otimes n} \, $  (see \S 2.1), thus  $ \, \delta_n(H) \subseteq
{J_{\scriptscriptstyle H}}^{\! \otimes n} = \h^n {\big( \h^{-1}
J_{\scriptscriptstyle H} \big)}^{\otimes n} \!\subseteq \h^n
{\big( H^\vee \big)}^{\otimes n} $,  which gives  $ \, H
\subseteq \! {\big( H^\vee \big)}' $.  For the rest, let
$ \, I' \! := \hbox{\sl Ker} \Big( H' {\buildrel \epsilon \over
{\relbar\joinrel\twoheadrightarrow}}\, R \,{\buildrel {\h \mapsto 0}
\over {\llongtwoheadrightarrow}}\, \Bbbk \Big) $;  \; as  $ \,
{\big(H'\big)}^{\!\vee} := \bigcup_{n=0}^\infty {\big( \h^{-1} I'
\big)}^n \, $,  \, to show that  $ \, H \supseteq {\big( H' \big)}^{\!
\vee} \, $  it is enough to prove that  $ \, H \!\supseteq \h^{-1} I'
\, $.  So let  $ \, x' \in I' \, $:  then  $ \, \delta_1(x') \in \h
\, H $,  \, hence  $ \, x' = \delta_1(x') + \epsilon(x') \in \h \,
H \, $,  \, thus  $ \, \h^{-1} x' \in H \, $.
                                                \par
   In the end, the last two identities follow directly from the
inclusions that we just proved.   \qed
\enddemo

\vskip9pt

\proclaim {Theorem 3.4} \, Let  $ \, H \in \HA \, $.
Then  $ \; H^\vee \in \QrUEA \, $.
\endproclaim

\demo{Proof}  By definition  $ H^\vee $  is generated by  $ \, J^\vee
:= \h^{-1} J \, $,  \, hence  $ H^\vee{\big|}_{\h=0} $  is generated
by  $ J^\vee{\big|}_{\h=0} $.  Now pick  $ \, j^\vee \in J^\vee \, $
and  $ \, j := \h \, j^\vee \in J \, $;  \, then  $ \; \Delta(j) =
\delta_2(j) + j \otimes 1 + 1 \otimes j - \epsilon(j) \cdot 1 \otimes
1 \in j \otimes 1 + 1 \otimes j + J \otimes J \, $,  \, therefore
$ \; \Delta \big( j^\vee \big)
%
%
 \, \in \, j^\vee \otimes 1 + 1 \otimes
j^\vee + \h^{+1} J^\vee \otimes J^\vee \, $,  \; whence  $ \;
\Delta \big(\overline{j^\vee}\big) = \overline{j^\vee} \otimes 1
+ 1 \otimes \overline{j^\vee} \; $  for  $ \, \overline{j^\vee}
:= j^\vee \mod \h \, H^\vee \in H^\vee{\big|}_{\h=0} \, $.  Thus
$ H^\vee{\big|}_{\h=0} $  is generated by its primitive part
$ P \big( H^\vee{\big|}_{\h=0} \big) $,  \, hence it is cocommutative.
Even more, this fact enables us to apply Lemma 5.5.1 in [Mo]
to  $ H^\vee{\big|}_{\h=0} \, $,  which then proves that
$ H^\vee{\big|}_{\h=0} $  is connected (as a Hopf algebra).
Thus  $ H^\vee{\big|}_{\h=0} $  is cocommutative, connected
and generated by its primitive part: that is, it is an rUEA,
so that  $ \, H^\vee \in \QrUEA \, $.   \qed
\enddemo

\vskip9pt

\proclaim {Theorem 3.5} \, Let  $ \; H \in \HA \, $.  Then
$ \; H' \in \QFA \, $.
\endproclaim

\demo{Proof} First,  $ \, H'{\big|}_{\h=0} \, $  is commutative as
a consequence of Lemma 3.2{\it (b)\/}  (cf.~[KT], Proposition 3.5).
                                          \par
   Second, we show that  $ \, {(H')}_\infty = {(I')}^\infty \, $.
For later use, set  $ \, I \! := I_{\scriptscriptstyle \! H} \, $,
$ J \! := \! J_{\scriptscriptstyle \! H} \, $,  $ J' \! := \!
J_{\scriptscriptstyle \! H'} \, $,  $ I' \! := \!
I_{\scriptscriptstyle \! H'} \, $.
                                            \par
   By definition  $ \, \h \, H' \! \subseteq \! I' \, $,  whence
$ \, H'_\infty := \bigcap_{n=0}^{+\infty} \h^n H' \! \subseteq
\! \bigcap_{n=0}^{+\infty} {\big( I' \big)}^n =: {\big( I'
\big)}^\infty $,  i.e.~$ {\big( H' \big)}_\infty \! \subseteq
\! {\big( I' \big)}^\infty $.  Conversely,  $ \, I' = \h \, H'
+ J' \, $  with  $ \, \h \, H' \subseteq \h \, H \, $  and  $ \,
J' = \delta_1(J') \subseteq \h \, H \, $:  \, thus  $ \, I' \!
\subseteq \! \h \, H \, $,  \, hence  $ \, {\big( I' \big)}^\infty
\! \subseteq \! \bigcap_{n=0}^{+\infty} \h^n H =: H_\infty \, $.
Now definitions give  $ \, H_\infty \! \subseteq \! H' \, $  and
$ \, \h^\ell H_\infty = H_\infty \, $  for all  $ \, \ell \in \Z
\, $,  \, so  $ \, \h^{-n} {\big( I' \big)}^\infty \! \subseteq \!
\h^{-n} H_\infty = H_\infty \! \subseteq \! H' \, $  hence  $ \,
{\big( I' \big)}^{\!\infty} \! \subseteq \! \h^n H' \, $  for all
$ \, n \in \N \, $,  thus finally  $ \, {\big( I' \big)}^{\!\infty}
\! \subseteq \! {\big( H' \big)}_\infty $.
                                            \par
   Third, we prove that  $ \, H'{\big|}_{\h=0} \, $  has no
non-trivial idempotents.  Let  $ \, a \in H' \, $  be such that
$ \; \overline{a} := a \, \mod \h \, H' \in H'{\big|}_{\h=0} \; $
is idempotent, i.e.~$ \, {\overline{a}}^{\,2} = \overline{a} \, $:
\, then  $ \, a^2 = a + \h \, c \, $  for some  $ \, c \in H' \, $.
Set  $ \, a_0 := \epsilon(a) \, $,  $ \, a_1 := \delta_1(a) \, $,
and  $ \, c_0 := \epsilon(c) \, $,  $ \, c_1 := \delta_1(c) \, $;
\, as  $ \, a $,  $ c \in H' \, $  we have  $ \, a_1 $,  $ c_1 \in
\h \, H \cap J = \h \, J \, $.  Applying  $ \delta_n $  to  $ \; a^2
= a + \h \, c \; $  we get  $ {\textstyle \sum\limits_{\Lambda \cup
Y = \{1,\dots,n\}}} \hskip-13pt \delta_\Lambda(a) \, \delta_Y(a) =
\delta_n\big(a^2\big) = \delta_n(a) + \h \, \delta_n(c) \; $  (for
all  $ \, n \in \N_+ \, $),  thanks to  Lemma 3.2{\it (a)\/}.  Since
$ \, a $,  $ c \in H' \, $  we have  $ \, \delta_n(a) $,  $ \delta_n(c)
\in \h^n H^{\otimes n} \, $  (for all  $ n \, $),  hence
 \vskip-17pt
  $$  \delta_n(a) \; \equiv  \hskip0pt  {\textstyle
\sum\limits_{\Lambda \cup Y = \{1,\dots,n\}}}  \hskip-9pt
\delta_\Lambda(a) \, \delta_Y(a)  \; = \;  2 \, \delta_0(a) \,
\delta_n(a) \; +  \hskip-5pt  {\textstyle \sum\limits_{\Sb
\Lambda \cup Y = \{1,\dots,n\} \\  \Lambda, Y \not= \emptyset
\endSb}}  \hskip-9pt  \delta_\Lambda(a) \, \delta_Y(a)
\hskip7pt \mod \, \h^{n+1} H^{\otimes n}  $$
 \vskip-5pt
for all  $ \, n \in \N_+ \, $,  \, which, recalling that  $ \, a_0
:= \delta_0(a) \, $,  \, gives (for all  $ \, n \in \N_+ $)
 \vskip-9pt
  $$  \big( 1 - 2 \, a_0 \big) \; \delta_n(a)  \hskip7pt  \equiv
\hskip7pt
   {\textstyle \sum_{\hskip-9pt  \Sb  \Lambda, Y \not= \emptyset  \\
\Lambda \cup Y = \{1,\dots,n\} \endSb}}  \hskip-1pt  \delta_\Lambda(a)
\, \delta_Y(a)  \hskip11pt  \mod \, \h^{n+1} H^{\otimes n} \; .
\eqno (3.1)  $$
 \vskip-3pt
   \indent   Now, applying  $ \, \epsilon \, $  to  $ \; a^2 = a + \h \,
c \; $  gives  $ \; {a_0}^{\!2} = a_0 + \h \, c_0 \, $.  This implies
$ \, \big( 1 - 2 \, a_0 \big) \not\in \h \, H \, $.  In fact, if  $ \,
\big( 1 - 2 \, a_0 \big) = \h \, \alpha \, $  ($ \, \alpha \in \h \,
H \, $),  \, then  $ \; \h^2 \, \alpha^2 = {(1 - 2 \, a_0)}^2 = 1
- 4 \, a_0 + 4 \, {a_0}^{\!2} = 1 + 4 \, \h \, c_0 \, $,  \; which
entails  $ \, 1 \in \h \, R \, $,  \, a contradiction.  From  $ \,
\big( 1 - 2 \, a_0 \big) \not\in \h \, H \, $  and (3.1), for all
$ \, n \in \N_+ \, $,  an easy induction gives  $ \, \delta_n(a)
\in \h^{n+1} H^{\otimes n} \, $  for all  $ n \, $.  Now take  $ \,
a_1 = a - a_0 = \h \, \alpha \, $  for some  $ \, \alpha \in \h \,
J \, $:  \, then  $ \, \delta_0(\alpha) = \epsilon(\alpha) = 0 \, $
and  $ \, \delta_n(\alpha) = \h^{-1} \delta_n(a) \in \h^n H^{\otimes n}
\, $  (for all  $ n \, $),  \, hence  $ \; \alpha \in H' \, $.  Thus
$ \; a = a_0 + \h \, \alpha \equiv a_0 \, \mod \h \, H' \, $,  \,
whence  $ \; \overline{a} = \overline{a_0} \in H'{\big|}_{\h=0}
\, $;  \, then  $ \; {\overline{a_0}}^{\,2} = \overline{a_0} \in
\Bbbk \, $  gives us  $ \, \overline{a_0} \in \{0,1\} \, $,  \,
hence  $ \, \overline{a} = \overline{a_0} \in \{0,1\} \, $,  \, q.e.d.
                                            \par
   Finally, assume that  $ \, p := \hbox{\it Char}\,(\Bbbk) > 0 \, $;
\, then we have to show that  $ \, \overline{\eta}^{\,p} = 0 \, $  for
each  $ \, \overline{\eta}^{\,p} \in J_{H'{|}_{\h=0}} \, $,  \, or
simply  $ \, \eta^p \in \h \, J_{\scriptscriptstyle H'} \, $  for
each  $ \, \eta \in J_{\scriptscriptstyle H'} \, $.  Indeed, for
any  $ \, n \in \N \, $  by the multiplicativity of  $ \Delta^n $
and  $ \, \Delta^n(\eta) = \sum_{\Lambda \subseteq \{1,\dots,n\}}
\delta_\Lambda(\eta) \, $  (cf.~\S 2.1) we get that  $ \,
\Delta^n(\eta^p) = \Big( \sum\nolimits_{\Lambda \subseteq
\{1,\dots,n\}} \delta_\Lambda(\eta) \Big)^p \, $  belongs to
$ \, \sum\limits_{\Lambda \subseteq \{1,\dots,n\}} \hskip-11pt
{\delta_\Lambda(\eta)}^p \, + \hskip-7pt \sum\limits_{\Sb e_1,
\dots, e_p < p \\   e_1 + \cdots e_p = p \endSb} \hskip-11pt
{p \choose e_1, \dots, e_p} \hskip-3pt \sum\limits_{\Lambda_k
\subseteq \{1,\dots,n\}, \, \forall \, k} \hskip-3pt
\prod_{k=1}^p {\delta_{\Lambda_k} (\eta)}^{e_k} + \,
\h \, \sum\limits_{k=0}^{n-1} \sum\limits_{\Sb \Psi
\subseteq \{1,\dots,n\}  \\  |\Psi|=k  \endSb} \hskip-11pt
j_\Psi \big( {J_{\!{}_{H'}}}^{\!\!\otimes k} \big) + \,
\h \, {J_{\!{}_{H'}}}^{\!\!\!\otimes n} \; $
%
%
%
%
%
%
because  $ \, \delta_\Lambda(\eta) \in j_\Lambda \Big(
{J_{\!{}_{H'}}}^{\!\!\otimes |\Lambda|} \Big) \, $  (for
all  $ \, \Lambda \subseteq \{1,\dots,n\} \, $)  and
$ \, \big[J_{\!{}_{H'}},J_{\!{}_{H'}}\big] \subseteq
\h \, J_{\!{}_{H'}} \, $.  Then
 \vskip-13pt
  $$  {\textstyle
   \delta^n(\eta^p) = {(\id_{\scriptscriptstyle H} \! -
\epsilon)}^{\otimes n} \big( \Delta^n(\eta^p) \big)
\, \in \,  {\delta_n(\eta)}^p \, + \hskip-11pt
\sum\limits_{\Sb e_1, \dots, e_p < p \\   e_1 + \cdots e_p = p \endSb}
\hskip-9pt  {p \choose {e_1, \dots, e_p}}  \hskip-3pt  \sum\limits_{\cup_k
\Lambda_k = \{1,\dots,n\}}  \hskip-13pt  \prod_{k=1}^p
{\delta_{\Lambda_k}(\eta)}^{e_k}  + \, \h \,
{J_{\!{}_{H'}}}^{\!\!\otimes n} \; .}  $$
 \vskip-7pt
\noindent
 Now,  $ \, {\delta^n(\eta)}^p \in {\big( \h^n H^{\otimes n} \big)}^p
\subseteq \h^{n+1} H^{\otimes n} \, $  since  $ \, \eta \in H' \, $;
\, similarly  $ \, \prod\nolimits_{k=1}^p {\delta_{\Lambda_k}
(\eta)}^{e_k} \in \h^{\sum_k |\Lambda_k| \, e_k} H^{\otimes n}
\subseteq \h^n H^{\otimes n} \, $  if  $ \, \bigcup_{k=1}^n
\Lambda_k = \{1,\dots,n\} \, $.  Moreover,  $ \, {p \choose
{e_1, \dots, e_p}} \, $  (with  $ \, e_1 $,  $ \dots $,  $ e_p
< p \, $)  is
%
%
zero in  $ \, \Bbbk = R \big/ \h \, R \, $,  i.e.~$ \,
{p \choose e_1, \dots, e_p} \in \h \, R \; $:  \, then  $ \;
\sum_{\Sb e_1, \dots, e_p < p \\   e_1 + \cdots e_p = p \endSb}
{p \choose e_1, \dots, e_p}  \sum_{\cup_k \Lambda_k = \{1,\dots,n\}}
\, \prod_{k=1}^p {\delta_{\Lambda_k}(\eta)}^{e_k} \in \h^{n+1}
H^{\otimes n} \, $.  Finally,  $ \, J_{\!{}_{H'}} \subseteq \h
\, J_{\!{}_H} \, $  implies  $ \, \h \, {J_{\!{}_{H'}}}^{\!\!
\otimes n} \subseteq \h^{n+1} H^{\otimes n} \, $.  Thus  $ \,
\delta_n(\eta^p) \in \h^{n+1} H^{\otimes n} \, $  for all
$ n \, $,  so  $ \, \eta \in \h \, H' \, $.   \qed
\enddemo

\vskip1,1truecm
 \eject

\centerline {\bf \S \; 4 \  Drinfeld's functors on
quantum groups }

\vskip10pt

\proclaim{Lemma 4.1} \, Let  $ \, F_\h \in \QFA \, $,  \,
with  $ \, F_\h{\big|}_{\h=0} \, $  reduced.  Let  $ \, I :=
I_{\scriptscriptstyle F_\h} \, $,  \, let  $ \widehat{F_\h} $
be the  $ I $--adic  completion of  $ F_\h \, $,  \,
$ \Rhat $  the  $ \h $--adic  completion of  $ R \, $,
\, and  $ \widehat{I{\phantom|\!}^n} $  the  $ I $--adic
closure of  $ I^n $  in  $ \widehat{F_\h} \, $,  \, for
all  $ \, n \in \N \, $.
                                        \par
   (a) \;  $ \, \widehat{F_\h} \, $  is isomorphic as
an  $ \Rhat $--module  to a formal power series algebra,
say  $ \, \Rhat \big[ \big[ {\{Y_b\}}_{b \in \Cal{S}} \big]
\big] \, $.
%
%
                                        \par
   (b) \; Use a section  $ \, \nu \, \colon \, \Bbbk \,
\lhook\joinrel\longrightarrow\, R \; $  of the quotient map
$ \, R \longtwoheadrightarrow R \big/ \h \, R =: \Bbbk \, $
to identify (set-theo\-retically)  $ \, \widehat{F_\h} \cong
\Rhat \big[\big[ {\{Y_b\}}_{b \in \Cal{S}} \big]\big] \, $
with  $ \, \nu(\Bbbk) \big[\big[ \{Y_0\} \cup {\{Y_b\}}_{b \in
\Cal{S}} \big]\big] \, $,  \, where  $ \, \h \cong Y_0 \, $.
Then  $ \, \widehat{\big( \widehat{I} \;{\big)}^n} \, $  and
$ \, \widehat{I{\phantom|\!}^n} \, $  coincide with the set of
all formal series of (least) degree  $ n $  in the  $ Y_i $'s
($ \, i \in \{0\} \! \cup \Cal{S} \, $),  for all  $ \, n
\in \N \, $.
                                        \par
%
   (c) \; Let  $ \, \mu : F_\h \loongrightarrow
\widehat{F_\h} \, $  be the natural map.  Then  $ \; \mu(F_\h)
\bigcap \widehat{I{\phantom|\!}^n} = \mu(I^n) \, $  for all
$ \, n \in \N \, $.
\endproclaim

\demo{Proof} Let  $ \, F[G] \equiv F_\h{\big|}_{\h=0} \! := F_\h \big/
\h \, F_\h \, $,  \, and let  $ \, \widehat{F[G]} = F[[G]] \, $  be
the  $ \germ $--adic  completion of  $ F[G] $,  with  $ \, \germ \! :=
\! \germ_e \! = \! \text{\sl Ker}\, \big( \epsilon_{\scriptscriptstyle
F[G]} \big) \, $  (see \S 2.1).  Then  $ \, I \! = \! \pi^{-1}(\germ)
\, $, \, where  $ \pi $  is the map  $ \, \pi \, \colon F_\h \!
\relbar\joinrel\twoheadrightarrow \! F_\h \big/ \h \, F_\h = F[G] \, $,
and so  $ \pi $  induces a continuous (specialization) epimorphism
$ \, \widehat{\pi}
\, \colon \, \widehat{F_\h} \longtwoheadrightarrow F[[G]] \, $.
%
%
                                                \par
   By construction  $ \widehat{F_\h} $  is a topological
$ \Rhat $--module  with  $ \, {\big( {\widehat{F_\h}} \,
\big)}_\infty = \{0\} \, $.  Set also  $ \, J := \text{\sl
Ker}\, (\epsilon_{\scriptscriptstyle F_\h}) \, $.
                                                \par
   {\it (a)} \,  Let  $ \, {\{y_b\}}_{b \in \Cal{S}} \, $  be a
$ \Bbbk $--basis  of  $ \, \germ \Big/ \germ^2 =: \gerg^\times \, $;
\, then  $ \, F[[G]] \cong \Bbbk \big[\big[ {\{Y_b\}}_{b \in \Cal{S}}
\big]\big] \, $  since  $ F[G] $  is reduced.  For each  $ \, b \in
\Cal{S} $,  \, pick a  $ \, j_b \in \pi^{-1}(y_b) \bigcap J \, $  and
fix a section  $ \, \nu \, \colon \, \Bbbk \lhook\joinrel\longrightarrow
R \, $  as in  {\it (b)}.  Now define a continuous morphism of
$ \Rhat $--modules  $ \, \Psi \, \colon \, \Rhat \big[\big[
{\{Y_b\}}_{b \in \Cal{S}} \big]\big] \!\loongrightarrow \widehat{F_\h}
\, $  mapping  $ \, Y^{\,\underline{e}} := \prod_{b \in \Cal{S}}
Y_b^{\,\underline{e}(b)} \, $  to  $ \, j^{\,\underline{e}} :=
\prod_{b \in \Cal{S}} j_b^{\,\underline{e}(b)} \, $  for all
$ \, \underline{e} \in \N^{\Cal{S}}_f := \big\{ \sigma \! \in \!
\N^{\Cal{S}} \,\big\vert\, \hbox{$ \sigma(b) = 0 $} \; \text{\sl
for almost all} \;\; b \in \Cal{S} \,\big\} \, $  (hereafter,
monomials like these are  {\sl ordered\/}  w.r.t.~any fixed
order of the set  $ \Cal{S} \, $).  Using  $ \nu $,  we can
set-theoretically identify  $ \, \Rhat \cong \nu(\Bbbk)[[Y_0]]
\, $  (with  $ \, \h \cong Y_0 \, $),
     \hbox{whence a bijection  $ \,
\nu(\Bbbk) \big[\big[ Y_0 \cup {\{Y_b\}}_{b \in \Cal{S}} \big]\big]
\cong \Rhat \big[\big[ {\{Y_b\}}_{b \in \Cal{S}} \big]\big] \, $
arises.}
                                          \par
%
%
   Is is easy to check that  $ \Psi $  is surjective.  To show
that it is injective too, look at graded rings, namely  $ \, G_\h
\big( \, \widehat{F_\h} \, \big) := \oplus_{n=0}^{+\infty} \Big(
\h^n \widehat{F_\h} \big/ \h^{n+1} \widehat{F_\h} \Big) \cong
\Bbbk \big[ \,\overline{\h} \,\big] \! \otimes_\Bbbk \! \Big(
\widehat{F_\h} \Big/ \h \, \widehat{F_\h} \Big) \cong \big( F[[G]]
\big)[Y_0] \, $.  In addition, there is an epimorphism of  {\sl
$ \Rhat $--modules}  $ \; \Psi \, \colon \, \Rhat \otimes_\Bbbk
\! F[[G]] \relbar\joinrel\relbar\joinrel\twoheadrightarrow \!
\widehat{F_\h} \; $  which induces an isomorphism of  {\sl graded\/
$ \Bbbk $--algebras}  $ \, G_{Y_0}\big(\Rhat \otimes_\Bbbk F[[G]]\big)
\cong G_\h\big(\, \widehat{F_\h} \,\big) \, $:  \, then (cf.~[Bo],
Ch.~III, \S 2.8)  $ \, \Psi \, $  is an isomorphism too.
                                                \par
   {\it (b)} \,  Since  $ \, \widehat{I} = \widehat{\pi}^{-1} \big(
\text{\sl Ker}\,(\epsilon_{\scriptscriptstyle F[[G]]}) \big) =
\text{\sl Ker}\,(\epsilon_{\scriptscriptstyle \widehat{F_\h}}) +
\h \, \widehat{F_\h} \, $,  \, each element of  $ \widehat{I} $  is
expressed, via  $ \Psi $,  by a series of degree at least 1; moreover,
for all  $ \, b $,  $ d \in \Cal{S} \, $  we have  $ \, j_b \, j_d -
j_d \, j_b = \h \, j_+ \, $  for some  $ \, j_+ \in \text{\sl Ker}\,
(\epsilon_{\scriptscriptstyle \widehat{F_\h}}) \, $.  Therefore
when multiplying  $ n $  factors from  $ \widehat{I} $  expressed
by  $ n $  series of positive degree, we can reorder the unordered
monomials in the  $ y_b $'s  occurring in the multiplication process
and eventually get a formal series   --- with  {\sl ordered\/}
monomials ---   of degree at least  $ n \, $.
                                                \par
   {\it (c)} \,  The analysis above shows that the natural map
$ \, \mu : F_\h \longrightarrow \widehat{F_\h} \, $  induces
$ \Bbbk $--module  isomorphisms  $ \, {\big( \widehat{I} \;\big)}^n
\Big/ {\big( \widehat{I} \;\big)}^{n+1} \cong  \widehat{\big(
\widehat{I} \;{\big)}^n} \Big/  \widehat{\big( \widehat{I} \;
{\big)}^{n+1}} \cong \widehat{I^n} \Big/ \widehat{I{\phantom{|}
\!}^{n+1}} \, $  ($ \, n \! \in \! \N \, $),  so  $ \, G_I\big(
F_\h \big) = G_{\widehat{I}}\big(\,\widehat{F_\h}\,\big) \, $,
\, these being the graded algebras associated to the  $ I $--adic
and the  $ \widehat{I} $--adic  filtration.  Moreover, the given
description of the  $ \widehat{I^n} $'s implies  $ \, G_{\widehat{I}}
\big(\,\widehat{F_\h}\,\big) := \bigoplus\limits_{n=0}^{+\infty}
\widehat{I{\phantom{|}\!}^n} \Big/ \widehat{I{\phantom{|}\!}^{n+1}}
\cong \Bbbk \big[ Y_0, {\{Y_b\}}_{b \in \Cal{S}} \big] \, $  as
$ \Bbbk $--modules,  and the same for  $ \, G_I\big(F_\h\big) \, $.
It follows that  $ \, F_\h \big/ I^n \cong \widehat{F_\h} \Big/
\widehat{I{\phantom{|}\!}^n} $,  \, whence  $ \, \mu(I^n) =
\widehat{I{\phantom{|}\!}^n} \bigcap \mu(F_\h) \, $,  \,
as claimed.   \qed
\enddemo

%
%
%

\vskip7pt

\proclaim{Lemma 4.2}  Let  $ \, F_\h \in \, \QFA \, $,  \, and assume
that  $ F_\h{\big|}_{\h=0} $  is reduced.  Then:
                               \par
   (a) \, if  $ \, \varphi \in F_\h \, $  and
$ \, \h^s \, \varphi \in {I_{\scriptscriptstyle F_\h}
\phantom{)}}^{\!\!\!\!\!n} \, $  ($ s, n \in \N $),  \,
then  $ \; \varphi \in {I_{\scriptscriptstyle F_\h}
\phantom{)}}^{\!\!\!\!\!n-s} $;
                               \par
   (b) \, if  $ \, y \in I_{\scriptscriptstyle F_\h}
\setminus {I_{\scriptscriptstyle F_\h} \phantom{)}}^{\!\!\!\!\!2} $,
\, then  $ \; \h^{-1} y \not\in \h \, {F_\h}^{\!\vee} \, $;
                               \par
   (c) \,  $ \; {\big(F_\h^{\,\vee}\big)}_\infty =
{\big( F_\h \big)}_\infty \; \big( = {I_{\scriptscriptstyle
F_\h}}^{\!\!\!\infty} \,\big) \, $.
%
%
\endproclaim

\demo{Proof}  {\it (a)} \, Set  $ \, I := I_{\scriptscriptstyle F_\h}
\, $.  Consider  $ \, I^\infty := \bigcap_{\,n=0}^{+\infty} I^n \, $
and the quotient Hopf algebra  $ \, \overline{F}_\h := F_\h \big/
I^\infty \, $:  \, then  $ \, \bar{I} := I_{\scriptscriptstyle
\overline{F}_\h} = I \big/ I^\infty \, $.  By  Lemma 3.1{\it
(a)},  $ \overline{F}_\h $  is again a QFA,  with  $ \,
\overline{F}_\h{\big|}_{\h=0} \! = F_\h{\big|}_{\h=0} \, $
and  $ \, \bar{I}^\infty := {I_{\scriptscriptstyle \overline{F}_\h}
\phantom{)}}^{\hskip-7pt \infty} = \{0\} \, $.  Now,  $ \; \phi
\in I^\ell \Longleftrightarrow \overline{\phi} \in \bar{I}^\ell
\; $  for all  $ \, \phi \in F_\h \, $,  $ \, \ell \in \N \, $,  \,
with  $ \, \overline{\phi} := \phi + I^\infty \in \overline{F}_\h
\, $.  So it is enough to prove the claim for  $ \overline{F}_\h
\, $,  \, hence we can assume  $ \, I^\infty = \{0\} \, $;  \,
then the natural map from  $ F_\h $  to its  $ I $--adic
completion  $ \widehat{F_\h} \, $  is injective.  By the
proof of Lemma 4.1 one has  $ \, \widehat{I^\ell} \bigcap
F_\h = I^\ell $,  \, for all  $ \, \ell \, $:  \, then,
thanks to  Lemma 4.1{\it (b)},  $ \, \varphi \in F_\h \, $
and  $ \, \h^s \varphi \in \! I^n \, $  imply  $ \, \varphi
\in \widehat{I^{n-s}} \bigcap F_\h = I^{n-s} \, $.
                                                \par
   {\it (b)} \, Let  $ \, y \in I_{\scriptscriptstyle F_\h} \setminus
{I_{\scriptscriptstyle F_\h}}^{\!\!2} \, $.  Assume  $ \, \h^{-1} y =
\h \, \eta \, $  for some  $ \, \eta \in {F_\h}^{\!\vee} \setminus
\{0\} \, $.  As  $ \, {F_\h}^{\!\vee} := \bigcup_{N \geq 0} \h^{-N}
{I_{\scriptscriptstyle F_\h}}^{\!\!n} \, $  we have  $ \, \eta =
\h^{-N} i_N \, $  for some  $ \, N \in \N_+ \, $,  $ \, i_N \in
{I_{\scriptscriptstyle F_\h}}^{\!\!N} \, $.  Then  $ \,
\h^{-1} y = \h^{1-N} i_N \, $,  \, so  $ \,
\h^{N-1} y = \h \, i_N \, $:  \, but the r.-h.-s.~belongs
to  $ {I_{\scriptscriptstyle F_\h}}^{\!\!{N+1}} $,  whilst
the l.-h.-s.~cannot belong to  $ {I_{\scriptscriptstyle
F_\h}}^{\!\!{N+1}} $,  due to  {\it (a)\/}  and  $ \, y
\not\in {I_{\scriptscriptstyle F_\h}}^{\!\!2} \, $.
                                                \par
   {\it (c)} \,  $ \, F_\h \subseteq {F_\h}^{\!\vee} $  implies
$ {\big( F_\h \big)}_\infty \! \subseteq \! {\big( {F_\h}^{\!\vee}
\big)}_\infty $.  Conversely, by definition  $ {\big( F_\h
\big)}_\infty $  is a two-sided ideal of  $ F_\h $  and
$ {F_\h}^{\!\vee} $,  and  $ \, {\overline{F}_\h}^{\!\vee} \equiv
{\left( F_\h \Big/ {\big( F_\h \big)}_\infty \right)}^{\!\vee}
= {F_\h}^{\!\vee} \Big/ {\big( F_\h \big)}_\infty \, $,  \, so
$ \, {\big( {F_\h}^{\!\vee} \big)}_\infty \!\mod\! {\big( F_\h
\big)}_\infty \subseteq \! {\left( {\overline{F}_\h}^{\!\vee}
\right)}_\infty $  with  $ \, \overline{F}_\h := F_\h \Big/ \!
{I_{\scriptscriptstyle F_\h}}^{\!\!\!\infty} = F_\h \Big/ {\big(
F_\h \big)}_\infty \, $  (a QFA, by  Lemma 3.1{\it (a)}).  So
it's enough to show  $ \, {\left( {\overline{F}_\h}^{\!\vee}
\right)}_\infty \!\! = \{0\} \, $.
                                                \par
   Let  $ \, \mu \, \colon \, F_\h \longrightarrow \widehat{F_\h}
\, $  be as above: it embeds  $ \overline{F}_\h $  into  $ \,
\widehat{F_\h} \, $,  \, and gives  $ \, {\overline{F}_\h}^{\!
\vee} \! \subseteq {\widehat{F_\h}}^{\,\vee} \! := \bigcup_{n
\geq 0} \h^{-n} \widehat{I^n} \, $,  \, so  $ \, {\left(
{\overline{F}_\h}^{\!\vee} \right)}_{\!\infty} \! \subseteq
\! {\left( {\widehat{F_\h}}^{\,\vee} \right)}_{\!\infty} \, $.
By Lemma 4.1,  $ {\widehat{F_\h}}^{\,\vee} $  is contained in the
$ R $--subalgebra  of  $ \, F(R) \otimes_R \widehat{F_\h} \, $
generated by  $ \, {\big\{ \h^{-1} j_b \big\}}_{b \in \Cal{S}}
\, $,  \, which is  {\sl polynomial\/};  then  $ \, {\left(
{\widehat{F_\h}}^{\,\vee} \right)}_\infty \!\! = \{0\} \, $,
\, so  $ \, {\left( {\overline{F}_\h}^{\!\vee} \right)}_\infty
\!\! = \{0\} \; $  too.   \qed
%
%
\enddemo

\vskip7pt

\proclaim {Proposition 4.3}  Let  $ \, \text{\it Char}\,(\Bbbk)
= 0 \, $.  Let  $ \, F_\h \in \QFA \, $.  Then  $ \, {\big(
{F_\h}^{\!\vee} \big)}' \! = F_\h \, $.
\endproclaim

\demo{Proof}  Proposition 3.3 gives  $ \, F_\h \subseteq \big(
{F_\h}^{\!\vee} \big)' $,  \, and we must prove the converse.
Let  $ \, \overline{F}_\h := F_\h \Big/ {(F_\h)}_\infty \, $;
\, then  $ \, {\big( {F_\h}^{\!\vee} \big)}_\infty \! = {(F_\h)}_\infty \, $,  by  Lemma 4.2{\it (c)},  so
$ \, {\big( \overline{F}_{\!\h} \big)}^{\!\vee} \! =
%
%
 {F_\h}^{\!\vee} \! \Big/ {\big( {F_\h}^{\!\vee}
\big)}_\infty \, $,  \, and  $ \, \Big(\! \big( \overline{F}_{\!\h}
\big)^{\!\vee} \Big)' \!\! =
%
%
 \big( {F_\h}^{\!\vee} \big)' \! \Big/ \!
{(F_\h)}_\infty \, $.  Thus, if the claim is true for  $ \,
\overline{F}_{\!\h} \, $  then  $ \, F_\h \Big/ {(F_\h)}_\infty =:
\overline{F}_\h = \Big(\! \big( \overline{F}_{\!\h} \big)^{\!\vee}
\Big)' = \big( {F_\h}^{\!\vee} \big)' \Big/ {(F_\h)}_\infty \, $,
\, hence  $ \, \big( {F_\h}^{\!\vee} \big)' \! = F_\h \, $.  So
a proof for  $ \overline{F}_{\!\h} $  is enough,
     \hbox{thus we assume
$ \, I^\infty \! = {(F_\h)}_\infty \!\! = \! {\big( {F_\h}^{\!\vee}
\big)}_\infty \!\! = \{0\} \, $,  for  $ I \! := \!
I_{\scriptscriptstyle F_\h} $.}
                                               \par
   Let  $ \, x' \in {\big( {F_\h}^{\!\vee} \big)}' \, $;  \,
since  $ \, {(F_\h)}_\infty = \{0\} \, $  there are  $ \, n \in
\N \, $  and  $ \, x^\vee \in {F_\h}^{\!\vee} \setminus \h \,
{F_\h}^{\!\vee} \, $  such that  $ \, x' = \h^n x^\vee \, $.
By Theorem 3.4,  $ {F_\h}^{\!\vee} $  is a QrUEA
   \hbox{with  $ \,
{F_\h}^{\!\vee}{\big|}_{\h=0} \! = U(\gerg) \, $,  $ \, \gerg =
I^\vee \! \big/ \big( \h \, {F_\h}^{\!\vee} \bigcap I^\vee \big)
\, $,  and  $ \, I^\vee := \h^{-1} I \, $.}
                                             \par
   Fix an ordered  $ \Bbbk $--basis  $ \, {\{b_\lambda\}}_{\lambda
\in \Lambda} \, $  of  $ \gerg \, $,  \, a subset  $ \,
{\big\{ x^\vee_\lambda \big\}}_{\lambda \in \Lambda} \, $  of
$ {I_{\scriptscriptstyle F_\h}}^{\hskip-5pt \vee} $  such that
$ \, x^\vee_\lambda  \hskip-1,4pt  \mod \h \, {F_\h}^{\!\vee} =
b_\lambda \, $  for all  $ \, \lambda \in \Lambda \, $,  \, and
set  $ \, x_\lambda = \h \, x^\vee_\lambda \in J $,  \, for all
$ \lambda \, $.  If  $ \, d := \partial(\bar{x}) \, $  is the
degree of  $ \, \bar{x} \, $  w.r.t.~the standard filtration
of  $ \, U(\gerg) \, $,  \, then (by [EK], Lemma 4.12, or [KT],
\S 3.8)  $ \, d:= \partial(\bar{x}) \leq n \, $.  So we can
write  $ \overline{x^\vee} $  as a polynomial  $ \, P \big(
{\{b_\lambda\}}_{\lambda \in \Lambda} \big) \, $  of degree
$ \, d \leq n \, $;  \, hence  $ \, x^\vee = P \big( {\big\{ x^\vee_\lambda \big\}}_{\lambda \in \Lambda} \big) + \h \, x^\vee_{[1]} \, $  for some  $ \, x^\vee_{[1]} \in {F_\h}^{\!\vee} $.  Now  $ \, x' = \h^n P \big( {\big\{ x^\vee_\lambda \big\}}_{\lambda \in \Lambda} \big) + \h^{n+1} x^\vee_{[1]} \, $,  \, with  $ \, \h^n
P \big( \big\{ x^\vee_\lambda \big\}_{\lambda \in \Lambda} \big) \in F_\h \, $  since  $ P $  has degree  $ \, d \leq n \, $;  \, as
$ \, F_\h \! \subseteq \!\! {\big( {F_\h}^{\!\vee} \big)}' $
     \hbox{(by Proposition 3.3),  $ \, x'_1 := x' \! - \! \h^n P \big(\! {\big\{ x^\vee_\lambda \big\}}_{\! \lambda \in \Lambda}
\big) \in \! {\big( {F_\h}^{\!\vee} \big)}' \, $  and  $ \, x'_1
\! = \h^{n+1} x^\vee_{[1]} \! = \h^{n_1} x^\vee_1 \, $}
 for some  $ \, n_1 \in \N $,  $ \, n_1 > n \, $,  \, and  $ \,
x^\vee_1 \in {F_\h}^{\!\vee} \setminus \h \, {F_\h}^{\!\vee} $.  We
repeat this construction with  $ x'_1 $  instead of  $ x' $,  $ \,
n_1 $ instead of  $ n $,  etc.: iterating, we get an increasing
sequence  $ {\big\{ n_s \big\}}_{s \in \N} $  and a sequence  $ \,
{\big\{ P_s \big( {\{ X_\lambda \}}_{\lambda \in \Lambda} \big)
\big\}}_{s \in \N} \, $  such that  $ \, x' = \sum_{s \in \N}
\h^{n_s} P_s \big( {\big\{ x^\vee_\lambda \big\}}_{\lambda \in
\Lambda} \big) \, $  and the degree of  $ P_s \big( {\{ X_\lambda
\}}_{\lambda \in \Lambda} \big) $  is at most  $ n_s \, $.  By
construction  $ \, \h^{n_s} P_s \big( {\big\{ x^\vee_\lambda
\big\}}_{\lambda \in \Lambda} \big) \! \in {I_{\scriptscriptstyle
F_\h}}^{\hskip-3pt n_s} \, $
       \hbox{for all  $ s \, $,  \, so  $ \,
\sum_{s \in \N} \h^{n_s} P_s \big( {\big\{ x^\vee_\lambda
\big\}}_{\lambda \in \Lambda} \big) \in \widehat{F_\h} \, $  (the}
  $ I_{\! \scriptscriptstyle F_\h} \! $--adic  completion
of  $ F_\h $),  $ \, x' = \sum_{s \in \N} \h^{n_s} P_s \big(
{\big\{ x^\vee_\lambda \big\}}_{\lambda \in \Lambda} \big) \, $
is an identity in  $ \widehat{F_\h} \, $,  \, and  $ \, x' \!
\in \! {\big( {F_\h}^{\!\vee} \big)}' \bigcap \widehat{F_\h}
\, $.  Now consider the specialization map  $ \, \pi \, \colon \,
F_\h \relbar\joinrel\twoheadrightarrow F_\h{\big|}_{\h=0} = F[G]
\, $  and the embedding  $ \, \mu \, \colon \, F_\h \hookrightarrow
\widehat{F_\h} \, $:  \; then  $ \pi $  extends to  $ \, \widehat{\pi}
\, \colon \, \widehat{F_\h} \longtwoheadrightarrow \widehat{F[G]} =
F[[G]] \, $,  \, and  $ \, \mu{\big|}_{\h=0} \, \colon \, F[G] =
F_\h{\big|}_{\h=0} \llongrightarrow \widehat{F_\h}{\big|}_{\h=0} =
F[[G]] \, $  is injective too.  Since  $ \, \hbox{\sl Ker}\,(\pi) =
\h \, F_\h \, $  and  $ \, \hbox{\sl Ker}\,(\widehat{\pi}) = \h \,
\widehat{F_\h} \, $,  \, this implies  $ \, F_\h \bigcap \, \h \,
\widehat{F_\h} = \h \, F_\h \, $,  \, whence we get  $ \, F_\h
\bigcap \h^\ell \, \widehat{F_\h} = \h^\ell F_\h \, $  for all
$ \, \ell \in \N \, $.  Getting back to  $ \, x' \in {\big(
{F_\h}^{\!\vee} \big)}' \bigcap \widehat{F_\h} \, $,  \, we
have  $ \, x' = \h^{-n} y \, $  for some  $ \, n \in \N \, $
and  $ \, y \in F_\h \, $;  \, therefore  $ \, y = \h^n x'
\in F_\h \bigcap \h^n \, \widehat{F_\h} = \h^n F_\h\, $,
\, so eventually  $ \, x' \in F_\h \, $.   \qed
\enddemo

\vskip7pt

\proclaim{Proposition 4.4} Let  $ \, H, K \! \in \HA \, $,
\, and  $ \; \langle \,\ , \,\ \rangle \, \colon H \times K
\loongrightarrow R \, $  be a Hopf pairing.  Then
                                      \hfill\break
   \indent   \hbox{(a)  $ H^\vee \! \subseteq \!
{\big( K' \big)}^\bullet \! $  and  $ \, K' \! \subseteq
\! {\big( H^\vee \big)}^\bullet \! $  (and viceversa), and
the above induces a pairing  $ \displaystyle{ H^\vee \!\!
\times \! K' \!\! \rightarrow \! R } $.}
%
%
%
                                      \hfill\break
   \indent   (b) \, If the pairing and its specialization at
$ \h = 0 $  are both perfect, and  $ \, K \! = \! {H}^\bullet \, $,
then  $ \, K' \! = \! {\big( H^\vee \big)}^\bullet \, $.
\endproclaim

%
%
\demo{Proof}  {\it (a)} \, Scalar extension and restriction give a
pairing between  $ H^\vee $  and  $ K' \, $:  \, we must prove it
is  $ R $--valued.  Let  $ \, I = I_{\scriptscriptstyle H} \, $.
Pick  $ \, c_1 $,  $ \dots $,  $ c_n \in I \, $,  $ \, y \in K' \, $:
\, then  $ \; \big\langle \prod_{i=1}^n c_i \, , \, y \big\rangle =
\big\langle \! \otimes_{i=1}^n c_i \, , \, \Delta^n(y) \big\rangle =
\sum_{\Psi \subseteq \{1,\dots,n\}} \big\langle \! \otimes_{i=1}^n c_i
\, , \, \delta_\Psi(y) \big\rangle = \sum_{\Psi \subseteq \{1, \dots,
n\}} \big\langle \! \otimes_{i \in \Psi} c_i \, , \, \delta_{|\Psi|}(y)
\big\rangle \cdot {\textstyle \prod_{j \not\in \Psi}} \langle c_j \, ,
1 \rangle \, \in \, \sum_{\Psi \subseteq \{1,\dots,n\}} \h^{n - |\Psi|}
R \cdot \h^{|\Psi|} R = \h^n R \; $.  The outcome is  $ \; \big\langle
I^n, K' \big\rangle \subseteq \h^n R \, $,  \, whence  $ \; \big\langle
\h^{-n} I^n, K' \big\rangle \subseteq R \, $,  \, for all  $ \, n \in
\N \, $;  \, then  $ \, H^\vee \! \subseteq \! {\big( K' \big)}^\bullet $
and  $ \, K' \subseteq \! {\big( H^\vee \big)}^\bullet \, $,  \, hence
$ \, H^\vee \! \times K' \! \longrightarrow F(R) \, $  takes values
into  $ R \, $,  \, q.e.d.
                                             \par
   {\it (b)} \, Let  $ \, \psi \in {\big( H^\vee \big)}^\bullet \, $:
then  $ \, \big\langle \h^{-s} I^s , \, \psi \big\rangle \in R \, $,
\, so  $ \, \big\langle I^s , \, \psi \big\rangle \in \h^s R $,  for
all  $ s \, $.  For  $ s \! = \! 0 $  we get  $ \, \big\langle H, \,
\psi \big\rangle \in R \, $,  \, thus  $ \, \psi \in H^\bullet \! =
K \, $,  \,
and so  $ \,
\delta_n(\psi) \in K^{\otimes n} \, $  for all  $ n \, $.
Now,  $ \, H^\bullet = K \, $  implies  $ \, {\big( H^{\otimes n}
\big)}^\bullet = K^{\otimes n} \, $  w.r.t.~the induced pairing
$ \; H^{\otimes n} \times K^{\otimes n} \loongrightarrow R \, $.
Moreover, since  $ \, \big\langle 1_{\scriptscriptstyle H} \, ,
J_{\scriptscriptstyle K} \big\rangle = 0 \, $  and  $ \, \delta_n
(\psi) \in {J_{\scriptscriptstyle K}}^{\!\!\otimes n} \, $,  \,
inverting the previous argument we find  $ \, \big\langle
H^{\otimes n}, \, \delta_n(\psi) \big\rangle = \big\langle
{I_{\scriptscriptstyle H}}^{\!\!\otimes n}, \, \delta_n(\psi)
\big\rangle \subseteq \h^n R \, $  (for all  $ n \, $): thus
$ \, \h^{-n} \delta_n(\psi) \in {\big( H^{\otimes n} \big)}^\bullet
= K^{\otimes n} \, $,  \, i.e~$ \, \delta_n(\psi) \in \h^n
K^{\otimes n} \, $  for all  $ n \, $,  whence
$ \, \psi \in K' \, $.   \qed
\enddemo

\vskip7pt

\proclaim{Proposition 4.5} \, Let  $ \, \text{\it Char}\,(\Bbbk) = 0
\, $.  Let  $ \, U_\h \in \QrUEA \, $.  Then  $ \; {\big( {U_\h}' \,
\big)}^{\!\vee} = \, U_\h \; $.
\endproclaim

\demo{Proof}  If the claim holds for  $ \, \overline{U}_\h :=
U_\h \Big/ {\big( U_\h \big)}_\infty \, $,  \, then  $ \,
%
%
\overline{U}_\h = \Big( \! \big( \, \overline{U}_\h \,\big)'
\Big){\phantom{\big|}}^{\hskip-4pt \vee} \! = {\big( {U_\h}'
\big)}^{\!\vee} \!\! \Big/ \! {\big( U_\h \big)}_\infty \, $  by
%
%
 Lemma 3.1{\it (a)},  thus  $ \, {\big( {U_\h}'
\, \big)}^{\!\vee} \! = U_\h \, $.  Therefore a proof
for  $ \overline{U}_\h $  is enough, and we may assume
$ \, {\big( U_\h \big)}_\infty \! = \{0\} \, $.
                                             \par
   Our plan is to mimic the proof of the same result for
quantum groups ``\`a la Drinfeld'' in [Ga1], Proposition
3.4 (now we drop the
       \hbox{hypothesis  $ \, \dim(\gerg) \! < \!\! +\infty \, $,
with  $ \, U(\gerg) = U_\h \big/ \h \, U_\h \, $,  by [Ga1], \S 3.9).}
                                             \par
   To simplify notation, set  $ H := U_\h \, $.  Let  $ \widehat{H} $
and  $ \Rhat $  be the  $ \h $--adic  completion  of  $ H $  and
$ R \, $.  Then  $ \widehat{H} $  is a
%
%
topological Hopf  $ \Rhat $--algebra,  whose coproduct
takes values into  $ \, \widehat{H} \,
\widehat{\otimes} \, \widehat{H} := H \, \widehat{\otimes} \, H \, $,
\, the  $ \h $--adic  completion of  $ \, H \otimes H \, $;  \,
clearly  $ H $  embeds into  $ \widehat{H} $  as a (topological)
Hopf  $ R $--subalgebra.  Set also  $ \; \widehat{H}' := \big\{\,
\eta \in \widehat{H} \,\big\vert\, \delta_n (\eta) \in \h^n
\widehat{H}^{\,\widehat{\otimes}\, n} \, \big\} \, $,  \, and
$ \; \big( \widehat{H}' \big)^{\!\times} \! := \bigcup_{n \geq 0}
\h^{-n} {I_{\scriptscriptstyle \! \widehat{H}'}}^{\hskip-3pt n}
\; \Big( \! \subseteq Q \big( \Rhat \, \big) \otimes_{\Rhat}
\widehat{H} \, \Big) \, $,  \, where  $ \, I_{\scriptscriptstyle
\! \widehat{H}'} := \hbox{\sl Ker}\,(\epsilon_{\scriptscriptstyle
\! \widehat{H}'}) + \h \cdot \widehat{H}' \, $.  Finally, let
$ \, \big( \widehat{H}' \big)^{\!\vee} \, $  be the  $ \h $--adic
completion of  $ \big( \widehat{H}' \big)^{\!\times} \, $.
                                        \par
   Now consider  $ \, \widehat{K} := {\widehat{H}}^* \equiv
{\hbox{\sl Hom}}_{\Rhat} \big( \widehat{H} \, , \, \Rhat \,\big)
\, $:  \, it is a complete topological Hopf  $ \Rhat $--algebra,
w.r.t.{} the weak topology, in perfect Hopf pairing with
$ \widehat{H} \, $.  We set  $ \, \widehat{K}^\times
:= \sum_{n \geq 0} \h^{-n} {J_{\scriptscriptstyle \!
\widehat{K}}}^{\hskip-3pt n} \; \Big( \! \subseteq Q
\big( \Rhat \, \big) \otimes_{\Rhat} \widehat{K} \, \Big) $,
\, where  $ \, J_{\scriptscriptstyle \! \widehat{K}} :=
\hbox{\sl Ker}\,(\epsilon_{\scriptscriptstyle \! \widehat{K}})
\, $,  we let  $ \widehat{K}^\vee $  be the  $ \h $--adic
completion of  $ \widehat{K}^\times $,  \, and we define  $ \,
{\big( \widehat{K}^\vee \big)}' \, $  in the obvious way.  With
the same arguments as for Proposition 4.3, one proves  $ \, {\big(
\widehat{K}^\vee \big)}' = \widehat{K} \, $.  Like in [Ga1], one shows
(as for Proposition 4.4)  $ \, \widehat{H}' = {\big( {\widehat{K}}^\vee
\big)}^{\!\bullet} \, $  and  $ \, \widehat{K}^\vee \subseteq {\big(
{\widehat{H}}' \big)}^{\!\bullet} \, $;  \, moreover  $ \,
\widehat{H} = \widehat{K}^* $,  \, hence  $ \, \widehat{K}^\vee
= {\big( {\widehat{H}}' \big)}^{\!\bullet} \, $.  Using this and
$ \, {\big( \widehat{K}^\vee \big)}' = \widehat{K} \, $  one
proves  $ \, \big( \widehat{H}' \big)^{\!\vee} = \widehat{H}
\, $  too  (see [Ga1] for details).
                                        \par
   Definitions imply  $ \, \widehat{H} \big/ \h^n \widehat{H} =
H \big/ \h^n H \, $,  \, thus  $ \, \h^n \widehat{H} \bigcap H = \h^n
H \, $,  \, for all  $ \, n \in \N \, $;  \, similarly  $ \, \h^n
\widehat{H}^{\widehat{\otimes} \ell} \bigcap H^{\otimes \ell} = \h^n
H^{\otimes \ell} $,  \, for all  $ \, n $,  $ \ell \in \N \, $,  \,
whence  $ \, \widehat{H}' \bigcap H = H' \, $.  The description of
$ \widehat{H}' $  in [Ga1], \S 3.5 (which holds for  $ \, \dim(\gerg)
= \infty \, $  too),  tells us that  $ \, \widehat{H}' \bigcap \h^n
H = I{\phantom{|}\!}^{\,n}_{\!\!\widehat{H}'} \, $,  \, and also
(acting like in the proof of Lemma 4.1), since  $ \,
I_{\widehat{H}'} \bigcap H' = I_{\widehat{H}'} \bigcap H = I_{H'} \, $,
\, that  $ \, I{\phantom{|}\!}^{\,n}_{\!\!\widehat{H}'} \bigcap H =
I{\phantom{|}\!}^{\,n}_{\!\!H'} \, $  for all  $ \, n \in \N \, $.
                                        \par
   Finally, take  $ \, \eta \in H \setminus \h \, H \, $.
We can show that there is an  $ \, \eta' \in I^{\,\partial
(\overline{\eta})}_{\widehat{H}'} \, $  (notation of the
proof of Lemma 4.2)  such that  $ \, \eta' = \h^{\partial
(\overline{\eta})} \eta + \eta'_+ \, $  for some  $ \,
\eta'_+ \in I^{\,\partial(\overline{\eta}) + 1}_{\widehat{H}'}
\, $,  \, just like in [Ga1], \S 3.5.  Roughly, we consider any basis
of  $ \, H{\big\vert}_{\h=0} = \widehat{H}{\big\vert}_{\h=0} \, $
containing  $ \overline{\eta} \, $,  we look at the dual basis inside
$ \, \widehat{K}{\big\vert}_{\h=0} \, $  and lift it to a topological
basis of  $ K $,  then rescale the latter (dividing out each element
by a proper power of  $ \h \, $)  to sort a topological basis of
$ K^\vee $:  the dual basis of  $ \widehat{H}' $  will contain an
element  $ \eta' $  as required.  Then  $ \, \h^{\partial
(\overline{\eta})} \eta = \eta' - \eta'_+ \in I^{\,\partial
(\overline{\eta})}_{\widehat{H}'} \bigcap H = I^{\,\partial
(\overline{\eta})} \, $,  \, by the previous analysis:
so  $ \, \eta
%
%
 \in \h^{-\partial
(\overline{\eta})} I^{\,\partial(\overline{\eta})} \subseteq
{\big( H' \big)}^{\!\vee} \, $.  Then  $ \, H \subseteq {\big(
H' \big)}^{\!\vee} \, $,  \, whereas the reverse inclusion
follows from Proposition 3.3.   \qed
\enddemo

\vskip6pt

\proclaim {Corollary 4.6} \, Let  $ \, \text{\it Char}\,
(\Bbbk) = 0 \, $.  Let  $ \, U_\h \in \QrUEA \, $.  Then
$ \; {\big( {U_\h}' \big)}_F = {(U_\h)}_F \; $.
\endproclaim

\demo{Proof}  Definitions give  $ \, {H^\vee}_{\!\!F} \! = \! H_F \, $
for all  $ H $,  so  $ \, {\big( {U_\h}' \big)}_F \!\! = \! {\Big( \!
\big( {U_\h}' \big)^{\!\vee} \Big)}_{\!F} \!\! = \! {(U_\h)}_F \, $
by Proposition 4.5.   \qed
\enddemo

\vskip7pt

\proclaim{Theorem 4.7} \, Let  $ \, F_\h[G] \in \QFA \, $
(notation of Remark 1.5) such that  $ F_\h[G]{\big|}_{\h=0} $
is reduced.  Then we have  $ \; {F_\h[G]}^\vee{\Big|}_{\h=0} := \,
{F_\h[G]}^\vee \! \Big/ \h \, {F_\h[G]}^\vee = \, U(\gerg^\times)
\; $  as co-Poisson Hopf algebras (see\/ \S 1.1).
\endproclaim

\demo{Proof}  Set for short  $ \, F_\h := F_\h[G] \, $,  $ \,
F_0 := F_\h \big/ \h \, F_\h = F[G] \, $,  \, and  $ \, {F_\h}^{\!
\vee} := {F_\h[G]}^\vee $,  $ \, {F_0}^{\!\vee} := {F_\h}^{\!\vee}
\big/ \h \, {F_\h}^{\!\vee} \, $.  By Theorem 3.4,  $ \; {F_0}^{\!\vee}
= \U(\gerk) \; $  for some  $ \gerk \, $  (since  $ \Bbbk $  is a field!), and we want to improve this result.
                                             \par
   Again,  $ \, {\big( F_\h \big)}_\infty = {\big( {F_\h}^{\!\vee}
\big)}_\infty \, $  by  Lemma 4.2{\it (c)\/};  then  $ \,
\overline{{F_\h}^{\!\vee}} := {F_\h}^{\!\vee} \Big/ {\big(
{F_\h}^{\!\vee} \big)}_\infty = {F_\h}^{\!\vee} \big/ {\big(
F_\h \big)}_\infty = {\big( \overline{F_\h} \, \big)}^{\!\vee} \, $
by  Lemma 3.1{\it (a)},  and so  $ \, {F_\h}^{\!\vee}{\big|}_{\h=0}
\! = \overline{{F_\h}^{\!\vee}}{\big|}_{\h=0} \! = {\big(
\overline{F_\h} \, \big)}^{\!\vee}{\big|}_{\h=0} \, $  by
Lemma 3.1{\it (b)}.  Thus it is enough to prove the claim
for  $ \, \overline{F_\h} \, $,  \, i.e.~we can assume  $ \,
{\big( F_\h \big)}_\infty = \{0\} \, $.  Let  $ \, I :=
I_{\scriptscriptstyle F_\h} \, $,  \, and let  $ \widehat{F_\h} $
be the  $ I $--adic  completion of  $ F_\h \, $;  \, as  $ \,
I^\infty = {\big( F_\h \big)}_\infty = \{0\} \, $,  \, the
natural map  $ \, F_\h \loongrightarrow \widehat{F}_\h \, $
is a monomorphism.
                                             \par
   Let  $ \, J := \text{\sl Ker}(\epsilon_{F_\h}) \, $  and  $ \,
J^\vee := \h^{-1} J \subset {F_\h}^{\!\vee} \, $.  Let  $ \,
{\{y_b\}}_{b \in \Cal{S}} \, $  be a  $ \Bbbk $--basis  of  $ \,
J_0 \big/ {J_0}^{\!2} \, $,  \, with  $ \, J_0 := \text{\sl Ker}\,
(\epsilon_{\scriptscriptstyle F[G]}) \, $,  \, and lift it to  $ \,
{\{j_b\}}_{b \in \Cal{S}} \subseteq J \, $.  Using notation of Lemma
4.1,  $ \, I^n \big/ I^{n+1} \cong \widehat{I^n} \big/ \widehat{I^{n+1}}
\, $  (for all  $ n \, $);  \, then  Lemma 4.1{\it (b)\/}  implies that
$ \, I^n \big/ I^{n+1} \, $  has  $ \Bbbk $--basis  $ \, \big\{\, \h^{e_0}
j^{\,\underline{e}} \mod I^{n+1} \;\big|\; e_0 \in \N, \underline{e} \in
\N^{\,\Cal{S}}_f \, , \; e_0 + |\underline{e}| = n \,\big\} \, $  with
$ \, |\underline{e}| := \sum_{b \in \Cal{S}} \underline{e}(b) \, $.  As
$ \, \h^{-n} I^{n+1} = \h \cdot \h^{-(n+1)} I^{n+1} \equiv 0 \mod \h \,
{F_\h}^{\!\vee} \, $,  \, we argue that  $ \, \h^{-n} I^n \! \mod \h
\, {F_\h}^{\!\vee} \, $  is  $ \Bbbk $--spanned  by  $ \, \big\{\,
\h^{-|\underline{e}|} j^{\,\underline{e}} \mod \h \, {F_\h}^{\!\vee}
\;\big|\; \underline{e} \in \N^{\,\Cal{S}}_f \, , \; |\underline{e}|
\leq n \,\big\} \, $:  \, we claim this is a basis of  $ \, \h^{-n} I^n
\mod \h \, {F_\h}^{\!\vee} \, $.  If not, we find a non-trivial linear
combination which is zero: multiplying by  $ \h^n $  yields  $ \,
\gamma_n \in I^n \setminus I^{n+1} \, $  with  $ \, \h^{-n} \gamma_n
\equiv 0 \mod \h \, {F_\h}^{\!\vee} \, $;  \, then  $ \, \h^{-n}
\gamma_n \in \h \cdot \h^{-\ell} I^\ell \, $  for some  $ \, \ell
\in \N \, $,  \, so  $ \, \h^\ell \gamma_n = \h^{1+n} I^\ell
\subseteq I^{1 + n + \ell} \, $:  \, then  Lemma 4.2{\it (a)}
yields  $ \, \gamma_n \in I^{n+1} $,  \, a contradiction.
                                             \par
   Now let  $ \, j_\beta^\vee := \h^{-1} j_\beta \, $  for all
$ \, \beta \in \Cal{S} \, $.  Since  $ \, j_\mu \, j_\nu - j_\nu
\, j_\mu \in \h \, J \, $,  \, for any  $ \, \mu, \nu \in \Cal{S}
\, $,  \, we have  $ \; j_\mu \, j_\nu - j_\nu \, j_\mu = \h
\sum_{\beta \in \Cal{S}} c_\beta \, j_\beta + \h^2 \gamma_1 +
\h \, \gamma_2 \; $  for some  $ \, c_\beta \in R \, $,  $ \gamma_1
\in J \, $  and  $ \, \gamma_2 \in J^2 $,  \, whence  $ \; \big[
j_\mu^\vee, j_\nu^\vee \big] := j_\mu^\vee \, j_\nu^\vee -
j_\nu^\vee \, \gamma_\mu^\vee
%
%
\equiv \sum_{\beta \in \Cal{S}} c_\beta \,
j_\beta^\vee \; $  mod  $ \, J + J^\vee J \, $:
\, but  $ \, J + J^\vee J = \h \, \big( J^\vee + J^\vee
J^\vee \big) \subseteq \h \, {F_\h}^{\!\vee} \, $,  \, so  $ \;
\big[ j_\mu^\vee, j_\nu^\vee \big] \equiv \sum_{\beta \in \Cal{S}}
c_\beta \, j_\beta^\vee \mod \h \, {F_\h}^{\!\vee} \, $,  \, hence
$ \, \gerh := J^\vee \mod \h \, {F_\h}^{\!\vee} \, $  is a Lie
subalgebra of  $ {F_0}^{\!\vee} $.  But the latter has
$ \Bbbk $--basis  $ \, \big\{ {\big(j^\vee\big)}^{\underline{e}}
\mod \h \, {F_\h}^{\!\vee} \;\big|\; \underline{e} \in \N^{\,
\Cal{S}}_f \,\big\} \, $,  \, hence the PBW theorem gives  $ \,
{F_0}^{\!\vee} = U(\gerh) \, $  as algebras.  Also, the proof of
Theorem 3.4 gives  $ \; \Delta \big( j^\vee \big) \equiv j^\vee
\otimes 1 + 1 \otimes j^\vee \mod \, \h \, {\big( {F_\h}^{\!\vee}
\big)}^{\otimes 2} \; $  for  $ \, j^\vee \in J^\vee $,  \, so
$ \; \Delta(\text{j}) = \text{j} \otimes 1 + 1 \otimes \text{j}
\; $  for  $ \, \text{j} \in \gerh \, $,  \, whence  $ \,
{F_0}^{\!\vee} = U(\gerh) \, $  as  {\sl Hopf\/}  algebras too.
                                             \par
   Now, the specialization  $ \; \pi^\vee \colon \, {F_\h}^{\!
\vee} \relbar\joinrel\twoheadrightarrow {F_0}^{\!\vee} =
U(\gerh) \; $  restricts to  $ \; \eta \, \colon \,
J^\vee \relbar\joinrel\twoheadrightarrow \; \gerh
\, := J^\vee \!\!\! \mod \h \, {F_\h}^{\!\vee} = J^\vee \big/
J^\vee \cap {\big( \h \, {F_\h}^{\!\vee} \big)} = J^\vee \big/
\big( J + J^\vee J_\h \big) \, $,  \, for  $ \, J^\vee \cap
{\big( \h \, {F_\h}^{\!\vee} \big)} = J^\vee \cap \h^{-1}
{I_{\scriptscriptstyle F_\h}}^{\hskip-3pt 2} = J_\h + J^\vee
J_\h \, $  by  Lemma 4.2{\it (b)\/}.  Let  $ \; \rho \, \colon
\, J_0 \relbar\joinrel\twoheadrightarrow J_0 \big/ {J_0}^{\!2}
=: \gerg^\times \, $  be the projection, and  $ \; \nu \,
\colon \, \gerg^\times \lhook\joinrel\longrightarrow J_0 \, $
a section of  $ \rho \, $.  The specialization  $ \; \pi \,
\colon \, F_\h \relbar\joinrel\twoheadrightarrow F_0 \; $
restricts to  $ \; \pi' \colon \, J \relbar\joinrel\twoheadrightarrow
J \big/ (J \cap \h \, F_\h) = J_\h \big/ \h \, J_\h = J_0 \, $:
\, we fix a section  $ \, \gamma \, \colon \, J_0
\lhook\joinrel\longrightarrow J_\h \, $  of  $ \pi' $.
Also, multiplication by  $ \h^{-1} $  yields an  $ R $--module
isomorphism  $ \, \mu \, \colon \, J \, {\buildrel \cong \over
{\lhook\joinrel\relbar\joinrel\twoheadrightarrow}} \, J^\vee $.
                                             \par
   The composition  $ \, \sigma := \eta \circ \mu \circ \gamma \circ
\nu \; \colon \, \gerg^\times \longrightarrow \gerh \, $  is a vector
space isomorphism, independent of the choice of  $ \nu $  and
$ \gamma \, $:  \, we show it is also a Lie bialgebra isomorphism.
Pull-back via  $ \sigma $  the Lie bialgebra structure of  $ \gerh $
onto  $ \gerg^\times $,  and denote it by  $ \big( \gerg^\times,
{[\,\ , \ ]}_\bullet, \delta_\bullet \big) $;  \, also, denote by
$ \big( \gerg^\times, {[\,\ ,\ ]}_\times, \delta_\times \big) $
the Lie bialgebra structure dual to that of  $ \gerg \, $:  \,
we shall prove that these two structures coincide.
                                             \par
   First,  {\it for all  $ \, x_1 $,  $ x_2 \in \gerg^\times \, $
we have}  $ \; {\big[x_1,x_2\big]}_\bullet = {\big[ x_1, x_2
\big]}_\times \, $.  Indeed, let  $ \, f_i := \nu(x_i) \, $,
$ \, \varphi_i := \gamma(f_i) \, $,  $ \, \varphi^\vee_i :=
\mu(\varphi_i) \, $,  $ \, y_i := \eta \big( \varphi^\vee_i \big)
\, $  \, ($ \, i = 1, 2 $).  Then  $ \; {\big[ x_1, x_2 \big]}_\bullet
:= \sigma^{-1} \big( {\big[ \sigma(x_1), \sigma(x_2) \big]}_\gerh \big)
= \sigma^{-1} \big( [y_1,y_2] \big)  = \big( \rho \circ \pi' \circ
\mu^{-1} \big) \big( \big[ \varphi^\vee_1, \varphi^\vee_2 \big] \big)
= \big( \rho \circ \pi' \big) \big( \h^{-1} [\varphi_1, \varphi_2]
\big) = \rho \big( \{f_1,f_2\} \big) =: {\big[ x_1, x_2 \big]}_\times
\; $.
                                                 \par
   The case of Lie cobrackets is similar, using maps  $ \, \nu_\otimes
:= \nu^{\otimes 2} $,  $ \, \gamma_\otimes := \gamma^{\otimes 2} $,
etc., and notation  $ \, \chi_\otimes \! := \eta_\otimes \! \circ
\mu_\otimes = {(\eta \circ \mu)}^{\otimes 2} \, $  and  $ \, \nabla
\! := \Delta - \Delta^{\text{op}} $.  Now,  {\it for all  $ \, x \in
\gerg^\times \, $  we have}  $ \; \delta_\bullet(x) = \delta_\times(x)
\, $.  Indeed, let  $ \, f := \nu(x) \, $,  $ \, \varphi := \gamma(f)
\, $,  $ \, \varphi^\vee := \mu(\varphi) \, $,  $ \, y := \eta
\big( \varphi^\vee \big) \, $.  Then  $ \; \delta_\bullet(x) \!
:= {\sigma_\otimes}^{\!-1} \big( \delta_\gerh (\sigma(x)) \big)
\! = {\sigma_\otimes}^{\!\!-1} \big( \delta_\gerh \big( \eta \big(
\varphi^\vee \big) \big) \big) \! = {\sigma_\otimes}^{\!\!-1} \big(
\eta_\otimes \big( \h^{-1} \nabla \big( \varphi^\vee \big) \! \big)
\! \big) \! = {\sigma_\otimes}^{\!\!-1} \big( \! {\big( \eta \circ
\mu \circ \gamma \big)}_\otimes \big( \nabla(f) \big) \! \big) \!
%
%
 = \! \rho_\otimes \big(
\nabla(\nu(x)) \big) \! = \delta_\times(x) \, $,  where the last
equality holds because  $ \delta_\times(x) $  is characterized in
$ \, \gerg^\times \otimes \gerg^\times \, $  by  $ \; \big\langle
u_1 \otimes u_2 \, , \, \delta_\times(x) \big\rangle = \big\langle
[u_1,u_2] \, , \, x \big\rangle \; $  for all  $ \, u_1, u_2 \! \in
\! \gerg \, $,  \, while  $ \, \big\langle [u_1, u_2] \, , \, x
\big\rangle = \big\langle [u_1,u_2] \, , \, \rho(f) \big\rangle =
\big\langle u_1 \otimes u_2 \, , \rho_\otimes \big( \nabla(\nu(x))
\big) \big\rangle \, . \quad \hskip-3pt \square $
\enddemo
%
%

\vskip7pt

\proclaim{Theorem 4.8}  Let  $ \, \text{\it Char}\,(\Bbbk) = 0 \, $.
Let  $ \, U_\h(\gerg) \in \QrUEA \, $  (notation of Remark 1.5).
Then we have  $ \; {U_\h(\gerg)}'{\Big|}_{\h=0} := \, {U_q(\gerg)}' \Big/ \h \, {U_\h(\gerg)}' \, = \, F\big[G^\star\big] \; $  as
Poisson Hopf algebras (see\/ \S 1.1).
\endproclaim

\demo{Proof}  By Theorem 3.5,  $ \, {U_\h(\gerg)}' \, $  is a
QFA,  with  $ \; {U_\h(\gerg)}' {\buildrel \,{\h \rightarrow 0}\,
\over \llongrightarrow}\, F[H] \, $;
%
%
%
\; we must show  $ \, H = G^\star \, $.
%
%
Theorem 4.7 applied to  $ {U_\h(\gerg)}' $  yields
$ {\big( {U_\h(\gerg)}' \big)}^{\!\vee} \! {\buildrel
{\h \rightarrow 0}\, \over {\relbar\joinrel\longrightarrow}}
U(\gerh^\times) \, $;  then  $ U(\gerg) {\buildrel \,{0 \leftarrow \h}\,
\over {\leftarrow\joinrel\relbar\joinrel\relbar}} U_\h(\gerg) =
{\big({U_\h(\gerg)}'\big)}^{\!\vee} \!{\buildrel \,
{\h \rightarrow 0}\, \over {\relbar\joinrel\longrightarrow}}
\, U(\gerh^\times) $
by Proposition 4.5, \, thus  $ \, \gerh^\times = \gerg \, $:  \,
therefore  $ \, \gerh := {\big( \gerh^\times \big)}^\star =
\gerg^\star \, $,  \, whence  $ \, H = G^\star \, $,  \, q.e.d.
\qed
\enddemo

\vskip7pt

\proclaim {Theorem 4.9} Let  $ \, \text{\it Char}\,(\Bbbk) = 0
\, $.  Let  $ \; \langle \,\ , \, \ \rangle \, \colon F_\h \times
U_\h \loongrightarrow R \, $  be a perfect Hopf pairing with
$ \, F_\h \! \in \! \QFA \, $,  $ \, U_\h \! \in \! \QrUEA \, $,
$ \, F_\h = {U_\h}^\bullet $,  $ \, U_\h = {F_\h}^\bullet $.  Then
$ \, {U_\h}' = {\big( {F_\h}^{\!\vee} \big)}^\bullet \, $  and
$ \, {F_\h}^{\!\vee} = {\big( {U_\h}' \big)}^\bullet \, $.
\endproclaim

\demo{Proof} The assumptions imply that the specialized Hopf pairing
$ \, F_\h{\big|}_{\h=0} \times U_\h{\big|}_{\h=0} \!\longrightarrow
\Bbbk \, $  is perfect as well: then by Proposition 4.4 we have only
to prove the inclusion  $ \, {F_\h}^{\!\vee} \supseteq {\big( {U_\h}'
\big)}^\bullet \, $.
                                              \par
   Let  $ \, \varphi \in {\big( {U_\h}' \big)}^\bullet \, $,  \,
chosen so that  $ \, \big\langle \varphi \, , {U_\h}' \big\rangle =
R \, $.  Since  $ \, {\big( {U_\h}' \big)}^\bullet \subseteq F(R)
\otimes_R F_\h = F(R) \otimes_R {F_\h}^{\!\vee} \, $,  \, there is
$ \, c \in R \setminus \{0\} \, $  such that  $ \, \varphi_+ := c
\, \varphi \in {F_\h}^{\!\vee} \setminus \h \, {F_\h}^{\!\vee} \, $:
\, it follows that  $ \, \big\langle \varphi_+ \, , {U_\h}'\big\rangle
= c \, R \, $.  If  $ \, F_\h = F_\h[G] \, $,  $ \, U_\h = U_\h(\gerg)
\, $,  then Theorems 4.7--8 give  $ \, {F_\h}^{\!\vee}{\big|}_{\h=0}
= U(\gerg^\times) \, $  and  $ \, {U_\h}'{\big|}_{\h=0} = F[G^\star]
\, $.  Thus there is  $ \, \overline{\eta} \in F[G^\star] \, $  such
that  $ \, \big\langle \varphi_+{\big|}_{\h=0} \, , \overline{\eta}
\, \big\rangle = 1 \, $,  \, hence there is  $ \, \eta \in {U_\h}'
\, $  (a lift of  $ \overline{\eta} \, $)  such that  $ \, \big\langle
\varphi_+ \, , \eta \big\rangle = 1 + \h \, \kappa \, $  for some
$ \, \kappa \in R \, $;  \, but  $ \, \big\langle \varphi_+ \, ,
\eta \big\rangle \in c \, R \, $,  \, thus  $ c \, $
divides  $ (1 + \h \, \kappa) $  in  $ R \, $.
                                              \par
   As  $ \, \varphi_+ \in {F_\h}^{\!\vee}
%
%
 \, $  we have  $ \, \varphi_+ = \h^{-n}
\varphi_0 \, $  for some  $ \, n \in \N \, $  and  $ \, \varphi_0 \in
I{\phantom{|}\!}^{\,n}_{\!\!F_\h} \, $;  \, therefore  $ \, \big\langle
\varphi_0 , {U_\h}' \big\rangle = c \, \h^n R \, $.  On the other hand,
as  $ \, U_\h = \big({U_\h}'\big)^\vee \, $  (Proposition 4.5)
each  $ \, y \in U_\h \, $  can be written as  $ \, y = \h^{-\ell}
y{}' \, $  for some  $ \, \ell \in \N \, $  and  $ \, y{}' \in {U_\h}'
\, $;  then  $ \, \big\langle \varphi_0 , y \big\rangle = c \,
\h^{n-\ell} \big\langle \varphi , y{}' \big\rangle \in R \bigcap
c \, \h^{n-\ell} R \, $  because  $ \, \big\langle \varphi_0 ,
y \big\rangle \in R \, $  and  $ \, \big\langle \varphi , y{}'
\big\rangle \in R \, $.  Now, if  $ \, \h^{n-\ell} \big\langle
\varphi , y{}' \big\rangle \not\in R \, $  then  $ \, n - \ell
< 0 \, $  and so  $ \h \, $  divides  $ c \, $.  Since  $ c \, $
divides  $ (1 + \h \, \kappa) $  we get an absurd, unless  $ c \, $
is invertible in  $ R \, $:  then  $ \, \varphi = c^{-1} \varphi_+
\in {F_\h}^{\!\vee} \, $.  Otherwise, we have always  $ \,
\h^{n-\ell} \big\langle \varphi , y{}' \big\rangle \in R \, $,
hence  $ \, \big\langle \varphi_0 , y \big\rangle \in c \, R
\, $  for all  $ \, y \in U_\h \, $;  \, thus  $ \, c^{-1} \varphi_0
\in {U_\h}^\bullet = F_\h \, $.  Let  $ \widehat{F_\h} $  be the
$ I_{F_\h} $--adic  completion of  $ F_\h \, $:  the natural map
$ \, \mu : F_\h \loongrightarrow \widehat{F_\h} \, $  has kernel
$ \, I{\phantom{|}\!}^\infty_{\!\!F_\h} = {(F_\h)}_\infty = \{0\} \, $
(for  $ \, F_\h \in \QFA \, $,  and  $ {(F_\h)}_\infty $  is contained
in the trivial left radical of the perfect pairing between  $ F_\h $
and  $ U_\h \, $;  so  $ \, F_\h \subseteq \widehat{F_\h} \, $.  Now
$ \, c^{-1} \varphi_0 \in F_\h \subseteq \widehat{F_\h} \, $  and
$ \, \varphi_0 \in I{\phantom{|}\!}^{\,n}_{\!\!F_\h} \subseteq
\widehat{I{\phantom{|}\!}^{\,n}_{\!\!F_\h}} \, $:  then
$ \, c^{-1} \varphi_0 \in \widehat{I{\phantom{|}\!}^{\,n}_{\!
\!F_\h}} \, $  by  Lemma 4.1{\it (a)--(b)},  hence  $ \, c^{-1}
\varphi_0 \in \widehat{I{\phantom{|}\!}^{\,n}_{\!\!F_\h}} \bigcap
F_\h = I{\phantom{|}\!}^{\,n}_{\!\!F_\h} \, $,  \, by  Lemma 4.1{\it
(c)}.  Thus  $ \, \varphi = c^{-1} \h^{-n} \varphi_0 \in h^{-n}
I{\phantom{|}\!}^{\,n}_{\!\!F_\h} \subseteq {F_\h}^{\!\vee}
\, $.   \qed
\enddemo

\vskip4pt

   At last, we can gather our partial results to prove the
main result, i.e.~Theorem 2.2:

\vskip7pt

\demo{$ \underline{\hbox{\it Proof of Theorem 2.2}} $}  Claim
{\it (a)\/}  is proved in \S\S 3.3--5,  {\it (b)\/}  follows from
Propositions 4.3 and 4.5,  {\it (c)\/}  holds by Theorems 4.7--8,
whereas Theorem 4.9 proves  {\it (d)}.  Finally, assume  $ \,
\Char(\Bbbk) = 0 \, $  and consider  $ \, \H \in \HA_F \, $.
If  $ \, H_{(f)} \! \in \! \QFA \, $
       \hbox{(for some prime  $ \, \h
\in \! R \setminus \{0\} \, $)  is an  $ R $--integer  form of
$ \H \, $,  then}
 $ H_{(f)}^\vee $  is an integer form too
(by the very definitions) and a QrUEA (at  $ \h \, $),  by
Proposition 3.3.  Conversely, if  $ \, H_{(u)} \! \in \! \QrUEA \, $
(for some prime  $ \, \h \in \! R \setminus \{0\} \, $)  is an
$ R $--integer  form of  $ \H \, $,  then  $ H_{(u)}' $  is
an integer form too (by Corollary 4.6) and a QFA (again at
$ \h \, $),  by Proposition 3.5.
       \hbox{This proves  {\it (e)}.   \qed}
\enddemo

\vskip7pt

{\bf 4.10 Final remarks:}
                                                 \par
   {\it (a)} \, {\sl The Global Quantum Duality Principle as
a ``Galois correspondence'' theorem.} \, Theorem 2.2 says that
Drinfeld's functors set mutually inverse Galois-like correspondences
from  $ \HA $  to itself.  When  $ (\h) $  is maximal, the subcategories of quantum groups then are those of fixed objects for the composition
of these correspondences.  Namely, the composed operator  $ \, {\big(
{(\ )}^\vee \big)}' = {(\ )}' \circ {(\ )}^\vee \, $  plays the role
of a ``closure operator'', and  $ \, {\big( {(\ )}' \big)}^\vee
= {(\ )}^\vee \circ {(\ )}' \, $  plays the role of a
``taking-the-interior operator'': then QFAs may be thought
   \hbox{of as ``closed objects''  and QrUEAs as ``open objects''
in  $ \HA \, $.}
                                                 \par
   {\it (b)} \, {\sl Duality between Drinfeld's functors}. \, For  $ \,
n \in \N \, $  let  $ \; \mu_n \, \colon \, {J_{\scriptscriptstyle H}}^{\!
\otimes n} \lhook\joinrel\longrightarrow H^{\otimes n} \,{\buildrel
{\; m^n} \over \loongrightarrow}\, H \, $  be the composition of the
embedding of  $ {J_{\scriptscriptstyle H}}^{\!\otimes n} $  into
$ H^{\otimes n} $  with the  $ n $--fold  multiplication (in  $ H \, $):
then  $ \mu_n $  is the ``Hopf dual'' to  $ \delta_n \, $.  As  $ \,
H^\vee := \sum_{n \in \N} \mu_n\big( \h^{-n} {J_{\scriptscriptstyle
H}}^{\!\otimes n}\big) \, $  and  $ \, H' := \bigcap_{n \in \N}
{\delta_n}^{\!-1}\big(\h^{+n} {J_{\scriptscriptstyle H}}^{\!
\otimes n} \big) \, $,  \, the two functors are built up as
``dual'' to each other (cf.~also part  {\it (d)\/}  of Theorem 2.2).
                                                 \par
   {\it (c)} \, {\sl Ambivalence \;
       \hbox{QrUEA  $ \leftrightarrow $ QFA  \; in  $ \HA_F \, $.}}
\; Part  {\it (e)} of Theorem 2.2 means that  {\sl some Hopf algebras
{\sl over\/}  $ F(R) $  might be thought of  {\it both\/}  as ``quantum function algebras''  {\it and\/}  as ``quantum enveloping algebras''\/}: examples are  $ U_F $  and  $ F_F $  for  $ \, U \in \QrUEA \, $  and
$ \, F \in \QFA \, $.
                                                 \par
   {\it (d)} \, {\sl Drinfeld's functors for algebras, coalgebras and
bialgebras}.  The definition of either of Drinfeld functors requires
only half of the notion of Hopf algebra.  In fact, one can define
$ (\ )^\vee $  for all ``augmented algebras'' (i.e., roughly
speaking, ``algebras with a counit'') and  $ (\ )' $  for all
``coaugmented coalgebras'' (roughly, ``coalgebras with a unit''),
and in particular for bialgebras: this yields nice functors, and
neat results extending the global quantum duality principle (cf.~[Ga2]).
                                                 \par
   {\it (e)} \, {\sl Generalizations}.  A good deal of the results
about Drinfeld's functors can be extended to the case when  $ \,
R \big/ \h R \, $  is  {\sl not\/}  a field, i.e.~the (prime) ideal
$ (\h) $  is not maximal.  Similar considerations hold for remark
{\it (d)\/}  too.

\vskip1,57truecm

\vskip21pt

\enddocument